\documentclass{jams-l}

\RequirePackage{enumitem}
\RequirePackage{graphicx}
\RequirePackage{wrapfig}
\RequirePackage{tikz}
\usetikzlibrary{calc,hobby}

\newtheorem{theorem}{Theorem}[section]
\newtheorem{lemma}[theorem]{Lemma}
\newtheorem{proposition}[theorem]{Proposition}

\theoremstyle{definition}
\newtheorem{definition}[theorem]{Definition}
\newtheorem{example}[theorem]{Example}

\theoremstyle{remark}
\newtheorem{remark}[theorem]{Remark}
\newtheorem{assumption}[theorem]{Assumption}

\numberwithin{equation}{section}

\providecommand{\abs}[1]{\lvert#1\rvert}
\providecommand{\norm}[1]{\lVert#1\rVert}

\begin{document}

\title[Localisation for constrained transports II]{Localisation for constrained transports II: applications}


\author{Krzysztof J. Ciosmak}
\address{Department of Mathematics, University of Toronto, Fields Institute for Research in Mathematical Sciences}
\curraddr{}
\email{k.ciosmak@utoronto.ca}
\thanks{Part of the research presented here was performed when the author was a post-doctoral fellow at the University of Oxford, supported by the European Research Council Starting Grant CURVATURE, grant agreement No. 802689. 
 Another part was performed when the author was a doctoral student at the University of Oxford, supported by St John’s College, Clarendon Fund and Engineering and Physical Sciences Research Council.
}


\date{}

\dedicatory{}

\begin{abstract}
We present a range of applications of localisation for constrained transports for pairs of probability measures in order with respect to a lattice cone.

These examples comprise irreducible convex paving for martingale transports in infinite-dimensional spaces, irreducible convex paving for submartingale transports, localisation of the Monge--Kantorovich problem, pavings for harmonic transports, pavings for solutions of martingale problems.

Furthermore, we consider approximating martingale transports and characterise the polar sets in this setting, identifying them as the trivial polar sets -- the sets whose projections have vanishing measures. This provides a negative answer to a question of Cox.
\end{abstract}

\maketitle
\tableofcontents
 \tableofcontents
\section{Introduction}\label{s:intro}

This paper aims to present, illustrate, and exemplify applications of the theory developed in the accompanying paper \cite{Ciosmak20242}. The theory focuses on extending the concept of the irreducible convex paving of the martingale transport theory, which was discovered in the one-dimensional setting by Beiglb\"ock, Nutz, and Touzi in \cite{Nutz2017} and developed in the finite-dimensional setting by De March and Touzi in \cite{Touzi2019}. These developments were inspired by the applications of the martingale optimal transport problem to mathematical finance, introduced in \cite{Beiglbock2013} by Beiglb\"cok, Henry-Labord\`ere and Penkner. Martingale transports are however also of independent mathematical interest, e.g. in the context of Skorkhod embedding problem \cite{Beiglbock2017, Kallblad2017, Stebegg2022}.
The results concerning the paving say that, roughly speaking, any martingale transport between two given probability measures $\mu,\nu$ is constrained by the finest collection of closed, convex sets, whose relative interia partition the underlying space up to a set of $\mu$-measure zero. That is, if we are to transport mass from the relative interior of an element of the partitioning, then we are only allowed to do that within the closure of the element. Moreover, within the elements of the paving, the mass can be transported without further such restrictions. 

Our theory extends these results to the setting of transports between two given probability measures $\mu,\nu$ on a set $\Omega$ that obey constraints imposed by an arbitrary convex cone $\mathcal{F}$ of functions on $\Omega$ that is stable under taking finite maxima.
Namely, given any two probabilities $\mu,\nu$ in $\mathcal{F}$-order, i.e. $\int_{\Omega}f\,d\mu\leq\int_{\Omega}f\,d\nu$ for all $f\in\mathcal{F}$, we show that there exists the finest partitioning of $\Omega$ into sets that are $\mathcal{F}$-convex such that any $\mathcal{F}$-transport between $\mu$ and $\nu$ is constrained within the closures of the elements of the partitioning. 
Moreover, we show that if we consider transportation within the elements of the partitioning, then there are no further such restrictions. The notions of $\mathcal{F}$-convexity and $\mathcal{F}$-transports shall be defined in due course. 

Thanks to the above-mentioned extension, we are able to show that the existence and properties of the irreducible convex paving carries over to the infinite-dimensional setting, as well as to the infinite-dimensional submartingale transports, generalising results of Nutz and Stebegg \cite{Nutz2018}. We also include the localisation of the Monge--Kantorovich problem,  which found a wide range of applications in geometric and functional inequalities, see e.g. work of Klartag \cite{Klartag2017} and work of Cavalletti and Mondino \cite{Cavalletti2017,Mondino2017}. 
A novelty that follows from our results is the paving for harmonic transports, related to the Skorokhod embedding problem. We provide an example showing that partitioning for harmonic transports can be genuinely different than the convex paving for the same pair of measures. Another novel example is the paving for solutions to martingale problems.

The scope of applicability of our results is not limited to the above-mentioned settings but is much wider. We showcase these examples to illustrate the breadth of the results and flexibility of the general abstract setting.


\subsection{General setting}

Let $\Omega$ be a set, and let $\mathcal{F}$ be the lattice cone of functions on $\Omega$ generated by a set of functions $\mathcal{G}$.\footnote{We refer the reader to Section \ref{s:prelim} for the necessary definitions.} As an example illustrating the theory, the reader may consider $\mathcal{G}$ to be the set of affine, non-decreasing functions, or to be the set of all affine functions. In the first instance, this will lead to consideration of submartingale transports, whereas in the second this will lead to martingale transports.
We equip $\Omega$ with the weak topology $\tau(\mathcal{G})$ generated  by $\mathcal{G}$. We denote by $\mathcal{H}$ the complete lattice cone generated by $\mathcal{G}$.

Let $p\in\mathcal{H}$ be a proper function,  bounded from below by a positive constant.
Let $\mu,\nu$ be two Radon probability measures on $\Omega$, with respect to which $p$ is integrable.
We shall say that $\mu$ and $\nu$ are in $\mathcal{F}$-order, and write $\mu\prec_{\mathcal{F}}\nu$, provided that for any  $f\in\mathcal{F}$
\begin{equation*}
\int_{\Omega}f\, d\mu\leq\int_{\Omega}f\, d\nu.
\end{equation*}
An $\mathcal{F}$-transports between $\mu$ and $\nu$ is a Radon probability measure $\pi$ on $\Omega\times\Omega$ such that the respective marginals of $\pi$ are $\mu$ and $\nu$ and for all $f\in\mathcal{F}$ integrable with respect to $\mu$ and $\nu$ and all measurable, bounded and non-negative $h$ on $\Omega$ we have
\begin{equation*}
\int_{\Omega}h(\omega_1)f(\omega_1)\, d\mu(\omega_1)\leq \int_{\Omega\times\Omega}h(\omega_1)f(\omega_2)\, d\pi(\omega_1,\omega_2).
\end{equation*}
Equivalently, $\pi$ is an $\mathcal{F}$-transport between $\mu$ and $\nu$ if and only if, it is a distribution a pair of random variables $(X,Y)$, with $X\sim \mu$ and $Y\sim \nu$, such that for all $f\in\mathcal{F}$, $(f(X),f(Y))$ is a one-step submartingale. The set of such measures we shall denote by $\Gamma_{\mathcal{F}}(\mu,\nu)$.

Before we recall the results of \cite{Ciosmak20242} which we shall illustrate, let us state the standing assumptions.

\begin{assumption}\label{as:ass}
Let $\mathcal{G}$ be a set of functions on a set $\Omega$ that contains constants.
Let $\mathcal{F}$ be the lattice cone generated by $\mathcal{G}$. We shall assume that there exists $p$ in the complete lattice cone $\mathcal{H}$ generated by $\mathcal{F}$, that is non-negative and proper with respect to $\tau(\mathcal{G})$. 
Moreover, we assume also that $\mathcal{G}$ separates points of $\Omega$, consists of functions of $p$-growth. 
We shall denote by $\mathcal{A}$ the linear span of $\mathcal{G}$.
We shall assume that $\Omega$ is separable in $\tau(\mathcal{G})$.
\end{assumption}

The first result, \cite[Theorem 1.4]{Ciosmak20242}, concerns the existence of $\mathcal{F}$-transports. It generalises the classical theorem of Strassen, see \cite{Strassen1965, Ciosmak20232}.

\begin{theorem}\label{thm:strassenINTRO}
Suppose that $\mu,\nu$ are Radon, probability measures on $\sigma(\tau(\mathcal{G}))$ such that 
\begin{equation*}
\int_{\Omega}p\, d\mu<\infty\text{ and }\int_{\Omega}p\, d\nu<\infty
\end{equation*}
and that
\begin{equation}\label{eqn:majintro}
\int_{\Omega}f\,d\mu\leq\int_{\Omega} f\,d\nu
\end{equation}
for all $f\in\mathcal{F}$. Then there exists an $\mathcal{F}$-transport  $\pi\in\Gamma_{\mathcal{F}}(\mu,\nu)$.
\end{theorem}

The second result, \cite[Theorem 1.5]{Ciosmak20242} shows, for a given pair $\mu\prec_{\mathcal{F}}\nu$, the existence of the finest partitioning of $\Omega$ into $\mathcal{A}$-convex sets $\mathrm{irc}_{\mathcal{A}}(\mu,\nu)(\cdot)$ that constrains all elements of $\Gamma_{\mathcal{F}}(\mu,\nu)$, under the assumption that all transports between $\mu$ and $\nu$ are \emph{local}, meaning that the set within which each point is allowed to be transported possesses certain finite-dimensionality property. 
These sets shall be called \emph{irreducible components} and the partitioning shall be called \emph{irreducible}  $\mathcal{F}$-\emph{convex paving}.

Below $\mathcal{G}$ is said to be  symmetric whenever $\mathcal{G}=-\mathcal{G}$. Moreover, $\Phi\colon \Omega\to\mathcal{A}^*$ is the assignment to any $\omega\in\Omega$ the corresponding evaluation functional, given  by the formula
\begin{equation*}
    \Phi(\omega)(a)=a(\omega)\text{  for all }a\in\mathcal{A},\omega\in\Omega.
\end{equation*}
Here $\mathcal{A}^*$ is the dual space to $\mathcal{A}$, equipped with the norm
\begin{equation*}
    \norm{a}_{\mathcal{A}}=\sup\frac{\abs{a}}{p}(\Omega)\text{ for }a\in\mathcal{A}.
\end{equation*}
For a subspace $\mathcal{B}\subset\mathcal{A}$ we write $R_{\mathcal{B}}\colon\mathcal{A}^*\to\mathcal{B}^*$ for the restriction map, that assigns to any functional on $\mathcal{A}$ its restriction to $\mathcal{B}$. 

We also define $\Phi_{\mathcal{B}}\colon\Omega\to\mathcal{B}^*$, by the formula $\Phi_{\mathcal{B}}(\omega)(b)=b(\omega)$ for $b\in\mathcal{B}$ and $\omega\in\Omega$.

\begin{theorem}\label{thm:partitionINTRO}
Suppose that $\mu\prec_{\mathcal{F}}\nu$. Suppose that any $\mathcal{F}$-transport between $\mu$ and $\nu$ is local.
Then there exists a Borel measurable set $B\subset\Omega$ with $\mu(B)=1$ such that whenever $\omega_1,\omega_2\in B$ then 
\begin{equation*}
\mathrm{irc}_{\mathcal{A}}(\mu,\nu)(\omega_1)\cap\mathrm{irc}_{\mathcal{A}}(\mu,\nu)(\omega_2)=\emptyset
\end{equation*}
or
\begin{equation*}
\mathrm{irc}_{\mathcal{A}}(\mu,\nu)(\omega_1)=\mathrm{irc}_{\mathcal{A}}(\mu,\nu)(\omega_2).
\end{equation*}
These sets are $\mathcal{A}$-convex.
Moreover, for any $\lambda\in\Lambda_{\mathcal{F}}(\mu,\nu)$
\begin{equation*}
    \mathrm{supp}\lambda(\omega,\cdot)\subset \Phi^{-1}(\mathrm{cl}\mathrm{irc}_{\mathcal{A}}(\mu,\nu)(\omega))\text{ for }\mu\text{-almost every }\omega\in\Omega.
\end{equation*}
The sets $\mathrm{irc}_{\mathcal{A}}(\mu,\nu)(\cdot)$ are the smallest convex sets that satisfy  the above condition.
Furthermore, if $\mathcal{B}=\mathrm{Cone}\mathcal{G}\cap(-\mathrm{Cone}\mathcal{G})$, then
\begin{equation*}
  \Phi_{\mathcal{B}}(\omega)\in R_{\mathcal{B}}\big(\mathrm{irc}_{\mathcal{A}}(\mu,\nu)(\omega)\big)\text{ for }\mu\text{-almost every }\omega\in\Omega.
\end{equation*}
In particular, if $\mathcal{G}$ is symmetric, then
\begin{equation*}
    \omega\in \Phi^{-1}(\mathrm{irc}_{\mathcal{A}}(\mu,\nu)(\omega))\text{ for }\mu\text{-almost every }\omega\in\Omega.
\end{equation*}
\end{theorem}

Above, the set $\Lambda_{\mathcal{F}}(\mu,\nu)$ is the set of all disintegrations of $\mathcal{F}$-transports $\pi$ between $\mu$ and $\nu$, with respect to the first co-ordinate. We refer the reader to Section \ref{s:prelim} for the definitions and to \cite[Section 8]{Ciosmak20242} for detailed proofs of related auxiliary results.

In the third result, \cite[Theorem 1.8]{Ciosmak20242}, we consider polar sets with respect to all $\mathcal{F}$-transports between $\mu$ and $\nu$. These are defined as the Borel sets $U\subset\Omega\times\Omega$ such that they are of zero measure with respect to every $\mathcal{F}$-transport between $\mu$ and  $\nu$. 

\begin{theorem}\label{thm:polarinto}
Suppose that $\mu\prec_{\mathcal{F}}\nu$. Suppose that any $\mathcal{F}$-transport between $\mu$  and $\nu$ is local.
Let $U\subset\Omega\times\Omega$ be a Borel set such that for $\mu $-almost every $\omega\in \Omega$
\begin{equation*}
\Phi(U_{\omega})\subset\mathrm{irc}_{\mathcal{A}}(\mu,\nu)(\omega).
\end{equation*}
Then $U$ is $\Gamma_{\mathcal{F}}(\mu,\nu)$-polar set if and only if there exist Borel sets $N_1,N_2\subset\Omega$ with
\begin{equation*}
\mu(N_1)=0,\nu(N_2)=0
\end{equation*}
and 
\begin{equation*}
U\subset (N_1\times\Omega)\cup (\Omega\times N_2).
\end{equation*} 
\end{theorem}
 The theorem provides an affirmative answer to generalisation of a conjecture of Ob\l\'oj and Siorpaes \cite[Conjecture 1.3.]{Obloj2017}.

A description of the intersections of the closures of the irreducible components is provided by \cite[Theorem 13.11]{Ciosmak20242}. 

\begin{theorem}\label{thm:finestrucfinintro}
Suppose that $\Omega$ is separable in $\tau(\mathcal{G})$. There exists a Borel measurable set $B\subset\Omega$ with $\mu(B)=1$ such that whenever $\omega_1,\omega_2\in B$  are such that $\mathrm{clConv}_H\Phi(\mathrm{supp}\lambda(\omega_1,\cdot))$ and $\mathrm{clConv}_H\Phi(\mathrm{supp}\lambda(\omega_2,\cdot))$ are finite-dimensional and if 
\begin{equation*}
a^*\in\mathrm{clConv}_H\Phi(\mathrm{supp}\lambda(\omega_1,\cdot))\cap \mathrm{clConv}_H\Phi(\mathrm{supp}\lambda(\omega_2,\cdot)),
\end{equation*}
is such that for $i=1,2$ there exists $\eta_i\in\mathcal{P}_{\xi(p)}(\Omega)$ and $C_i>0$ such that $\eta_i\leq C_i\lambda(\omega_i,\cdot)$ and
\begin{equation*}
a^*(a)=\int_{\Omega}a\, d\eta_i\text{ for all }a\in\mathcal{A}.
\end{equation*}
then the Gleason parts of $a^*$ in
\begin{equation*}
\mathrm{clConv}_H\Phi(\mathrm{supp}\lambda(\omega_i,\cdot))
\end{equation*}
for $i=1,2$ are equal.
\end{theorem}

Let us also mention that \cite[Theorem 13.16]{Ciosmak20242} shows that the obtained partitioning is regular enough that measure $\mu$ may be disintegrated with respect to the partitioning, in such a way that almost every conditional measure has a single irreducible component.

\subsection{Related works}

We shall now describe in a greater detail the examples that our theory generalises.

\subsubsection{Martingale transport in finite-dimensional spaces}

 If   $\mathcal{G}$ is the linear space of all affine functions,  $\mathcal{H}$ is the cone of all convex functions on a finite-dimensional linear space, and $p$ is a norm on the space, then our setting corresponds to the martingale transport in finite-dimensional spaces. Works of Ghoussoub, Kim, Lim \cite{Ghoussoub2019}, De March, Touzi \cite{Touzi2019} and Ob\l\'oj and Siorpaes \cite{Obloj2017} are concerned with the convex paving for martingale transports.

The martingale transport problem in the one-dimensional discrete setting has been conceived by Henry-Labord\`ere, Beiglb\"ock and Penkner \cite{Beiglbock2013}, where it served  as the dual of the problem of robust, model-free, superhedging of exotic derivatives in financial mathematics. It has been further developed by Beiglb\"{o}ck, Juillet \cite{Juillet2016}, who studied the supports of martingale optimal transports, in a similar way  to the work of Gangbo and McCann \cite{McCann1996}.
In  \cite{Nutz2017} Beiglb\"ock, Nutz and Touzi introduce a quasi-sure formulation of the dual problem. It allowed the authors to obtain general duality result, with no duality gap, and prove the existence of dual optimisers. 
Further literature concerning the martingale transport problem includes  \cite{Galichon2014, DeMarch2018, DeMarch20182, Ekren2018, Dolinsky2015, Cheridito2021}. 

The duality in the martingale optimal transport problem remains an active area of research. In \cite{Backhoff2020} the authors develop a displacement interpolation for the martingale optimal transport, which is an analogue of the McCann's interpolation \cite{McCann1995, McCann1997} and the optimal Brenier map \cite{Brenier1991}. Further studies \cite{Backhoff2023, Backhoff20232} in this direction investigate displacement convexity of Bass functional, along generalised geodesics. 

In  Secion \ref{s:martingale} we show how our general theorems apply in  the  context of martingale transports on infinite-dimensional linear spaces.

\subsubsection{Supermartingale transport in one-dimensional spaces}\label{s:supmartingale}

Another example concerns  supermartingale transport. Let $\mathcal{G}$ be the space of all affine, non-decreasing functions  on the real line. Then $\mathcal{H}$ is the cone of  convex and non-increasing functions.
Two Borel probability measures $\mu,\nu$ on the real line with finite first moments are in $\mathcal{F}$-order if and only if they are distributions of the respective marginals of a one-step supermartingale, i.e., of a random variable $(X,Y)$ such that
\begin{equation*}
X\geq \mathbb{E}(Y\vert\sigma(X)).
\end{equation*}
In the one-dimensional case, this problem has been studied by Nutz and Stebegg in \cite{Nutz2018}. In \cite[Proposition 3.4.]{Nutz2018}, it is proven that given two Borel probability measures $\mu\prec_{\mathcal{F}}\nu$ on $\mathbb{R}$, with finite first moments, there exists a partitioning of $\mathbb{R}$ into open intervals such that any $\mathcal{F}$-transport between $\mu$ and $\nu$ has to occur within their closures. Moreover, a characterisation of polar sets for all $\mathcal{F}$-transports has been provided. 

We provide a detailed description of general submartigale transports in Section \ref{s:supermartingale}.

\subsubsection{Optimal transport localisation}

One of the aims of this paper is to illustrate the general localisation framework introduced in \cite{Ciosmak20242} and to show that it comprises both the irreducible convex paving pertaining to martingale transports and the localisation of the optimal transport theory.

The localisation in the context of Riemannian geometry has been developed by Klartag \cite{Klartag2017}, in the context of Finsler geometry by Ohta \cite{Ohta2018}, and by Cavalletti and Mondino \cite{Mondino2017, Cavalletti2017} to essentially non-branching metric measure spaces that satisfy curvature-dimension conditions.

Let us describe the method in the context of Riemannian geometry. Let $\mathcal{M}$ be a Riemannian manifold, equipped with a metric $d$ and the volume measure $\lambda$.
Let $\mu,\nu\in\mathcal{P}(\mathcal{M})$ be two Radon probability measures on $\mathcal{M}$, absolutely continuous with respect to $\lambda$ and with finite first moments.
Consider the Monge--Kantorovich transport problem between $\mu$ and $\nu$.  We refer the reader to  \cite{Monge1781, Kantorovich1958} for foundational works,  the books of Villani \cite{Villani2003, Villani2009}, lecture notes of Ambrosio \cite{Ambrosio2003}, and the book of Santambrogio \cite{Santambrogio2015} for a more recent account on optimal transport problems. 
The task is to find optimal coupling $\pi$ of $\mu$ and $\nu$ that minimises the total cost
\begin{equation*}
    \int_{\mathcal{M}\times\mathcal{M}}d(\omega_1,\omega_2)\, d\pi(\omega_1,\omega_2).
\end{equation*}
The Kantorovich duality tells that the minimal value of the above integral is equal to the supremum of all intergrals
\begin{equation*}
    \int_{\mathcal{M}}u\, d(\nu-\mu)
\end{equation*}
where $u\colon\mathcal{M}\to\mathbb{R}$ is a $1$-Lipschitz function. It can be shown that the supremum is attained.
Let $v$ be a $1$-Lipschitz function on $\mathcal{M}$ that is a Kantorovich potential, i.e., it attains the supremum. 
Let $\mathcal{F}$ be the tangent cone at $v$ to the set of $1$-Lipschitz functions on $\mathcal{M}$. It is easy to verify that $\mu\prec_{\mathcal{F}}\nu$. 

We say  that a set $\mathcal{T}\subset\mathcal{M}$ is a transport ray whenever it is a maximal set, with respect to inclusion, such that $v$ is isometric on $\mathcal{T}$. The transport rays partition $\mathcal{M}$, up to a set of $\lambda$-measure zero, see e.g. \cite{Klartag2017}. A crucial observation was that this partitioning inherits curvature-dimension properties of the original space and moreover, the masses of the conditional measures of $\mu$ and $\nu$ on the transport rays are equal. That allowed to provide proofs of e.g. L\'evy--Gromov isoperimetric inequality in a great generality. 

Let us note that it has been conjectured in \cite{Klartag2017} that the method can be generalised to multi-dimensional setting employing optimal transport of vector measures. The optimal transport of vector measures has been introduced in \cite{Ciosmak20211}, where one part of  this conjecture has been refuted. Another version of optimal transport of vector measures has been proposed in \cite{Ciosmak20212}. The other part of the conjecture of Klartag has been partially resolved in the affirmative in \cite{Ciosmak2021}. The method that was proposed in \cite{Klartag2017} for a resolution has been shown in \cite{Ciosmak20213} to contain a gap.

We shall show in Section \ref{s:optimal}, Proposition \ref{pro:subspace}, that the closure of the cone $\mathcal{F}$, when viewed as a subset of the space of integrable functions with respect to $\mu+\nu$, contains a linear subspace $\mathcal{A}$. As $\mu\prec_{\mathcal{F}}\nu$, we also see that $\mu\prec_{\mathcal{A}}\nu$. The subspace $\mathcal{A}$ consists of all functions that are constant on transport rays of $v$. 

This observation shows that there is a partitioning into $\mathcal{A}$-convex sets, which constrains all optimal transports between $\mu$ and $\nu$. In Section \ref{s:optimal} we provide a detailed description of this observation. 

\subsection{Generalised convexity}\label{s:convex}

 The elements of the irreducible partitioning exemplified in this paper and introduced in \cite{Ciosmak20242} are $\mathcal{A}$-convex sets.

If $\mathcal{F}$ is a convex cone of functions on a topological space $\Omega$, then a  closed set $K\subset\Omega$ is said to be $\mathcal{F}$-convex  whenever it equal to its $\mathcal{F}$-convex hull, that  is
\begin{equation*}
\mathrm{clConv}_{\mathcal{F}}K=\{\omega\in\Omega\mid f(\omega)\leq \sup f(K)\text{ for all }f\in\mathcal{F}\}.
\end{equation*}
Let us observe that when $\mathcal{F}$ is the complete lattice cone, or the lattice cone, generated by a linear space $\mathcal{A}$ then $K$ is $\mathcal{F}$-convex if and only if it is $\mathcal{A}$-convex. Moreover, if $\Omega$ is a locally convex topological vector space and $\mathcal{A}$ is the space of all continuous linear functionals on $\Omega$,  then the above defined notion of convexity coincides with the usual notion, as follows from the Hahn--Banach theorem.

The notion of generalised convexity with respect to an arbitrary family of functions was introduced by Fan in  \cite{Fan1963}. It allowed for a generalisation of the Krein--Milman theorem and has been developed and further generalised, see \cite{Dolecki1978, Rubinov2000, Singer1997}. These developments mainly were concerned with various dualities,  notions of subdifferentials and Fenchel transforms. These have found its application in the monopolist's problem by Figalli, Kim  and McCann \cite{Figalli2011} and by McCann and Zhang in \cite{McCann2019}.

We refer the reader to the book of H\"ormander \cite{Hormander2007} for an account comprising matter on notions of convexity with respect to: subharmonic functions, \cite[Chapter III]{Hormander2007}, plurisubharmonic  functions, \cite[Chapter  IV]{Hormander2007}. We note  also that the latter is equivalent to the notion of  holomorphic convexity of complex analysis, see e.g. \cite[Chapter 3]{Krantz2001}.

In \cite{Ciosmak2023} we study the generalised convexity in a way that is relevant to the current setting. The paper is concerned with a generalisation of the Levi problem \cite{Levi1911} and the Cartan--Thullen theorem \cite{Cartan1932}. In particular, theorems in \cite{Ciosmak2023} characterise topological spaces $\Omega$, equipped with a linear space of continuous functions on $\Omega$, such that there exists a proper function $p$ that belongs to the complete lattice cone $\mathcal{H}$ generated by $\mathcal{A}$. As it turns out, these spaces are precisely spaces that are complete with respect to $\mathcal{A}$. 

The notion of $\mathcal{F}$-convex sets is discussed in detail in Section  \ref{s:fconvex}.

\subsection{Balayages of measures}\label{s:balayage}

Let $\mathcal{F}$ be a lattice cone of Borel measurable functions on a topological space $\Omega$. For two Radon probability measures $\mu,\nu$ we write $\mu\prec_{\mathcal{F}}\nu$ whenever for any $f\in\mathcal{F}$ 
\begin{equation*}
\int_{\Omega}f\,d\mu\leq\int_{\Omega}f\, d\nu,
\end{equation*}
provided that the integrals are well-defined.

The above order on measures, related to the notion of balayage, originates in the works on potential theory of Poincar\'e. It was studied by Mokobodzki in \cite{Mokobodzki1984}, where he defines $\nu$ to be balayage of $\mu$ relative to $\mathcal{F}$ if $\mu\prec_{\mathcal{F}}\nu$.
We refer the reader also to \cite{Bowles2019, Ghoussoub2023, Ciosmak20232} for recent developments and investigations of balayage in the context of optimal transport. We refer also to \cite[Chapter XI]{Meyer1966} for an accessible introduction to the theory of balayages and for a proof of Strassen's theorem \cite[T 51, p. 245]{Meyer1966}.

We note that the balayage ordering  is a special case of stochastic ordering, studied e.g. in  \cite{Shaked2007}. The paper \cite{Kim2024} is related both  to the stochastic orderings and to the balayages of measures.

Examples of balayage ordering we shall provide along the presentation of our results.

\subsection{Novel results}

Let us give  an account of the new results that we obtain as applications of Theorem  \ref{thm:strassenINTRO}, Theorem \ref{thm:partitionINTRO}, and Theorem \ref{thm:polarinto}.

\subsubsection{Martingale transports in infinite-dimensional spaces}

The basic example is the martingale transport. Our theory allows for a treatment of the infinite-dimensional case, which generalises results  of \cite{Juillet2016, Nutz2017, Touzi2019, Obloj2017}. Specifically, suppose that $X^*$ is a dual, separable Banach space. Let $\mathcal{A}$ denote the space of all affine weakly* continuous functionals on $X^*$ and $\mathcal{H}$ be the cone of convex and lower semi-continuous, with respect to the weak* topology on $X^*$, functions on  $X^*$. We shall show in Section \ref{s:martingale}, Theorem \ref{thm:martingale}, that given 
two Radon probability measures $\mu,\nu$ on $X^*$ with finite first moments such that $\mu\prec_{\mathcal{H}}\nu$, and such that any $\mathcal{H}$-transport between $\mu$ and $\nu$ is local, see Definition \ref{def:local}, there exists the finest partitioning of $X^*$ into relatively open, convex and finite-dimensional subsets of $X^*$ that constrain all $\mathcal{H}$-transports between $\mu$ and $\nu$. We also provide a characterisation of polar sets, whose sections are contained in the corresponding irreducible components. This provides a positive resolution of a conjecture of Ob\l\'oj and Siorpaes posed in \cite[Conjecture 1.3.]{Obloj2017}.
In Theorem \ref{thm:intersections} we describe the fine structure of the intersections of the closures of the irreducible components.

\subsubsection{Submartingale transports in infinite-dimensional spaces}

We shall consider submartingale transport of measures on infinite-dimensional spaces, generalising the results of Nutz and Stebegg \cite{Nutz2018} on supermartingale transports described in Section \ref{s:supmartingale}.




We extend these results to the infinite-dimensional setting. Theorem \ref{thm:submartingale} shows the existence of the finest partitioning and provides the characterisation of the polar sets, while Theorem \ref{thm:intersectionssub} provides a description of the fine structure of the intersections of the closures of the irreducible components. Our treatment includes the case of martingale-submartingale transports, i.e., distributions that are martingales in certain directions and submartingales in the other directions.


The submartingale transport is relevant to the studies of the monopolist's problem investigated  by Rochet and Chon\'e \cite{Rochet1998} and, more recently, by Figalli, Kim, McCann \cite{Figalli2011} and McCann, Zhang \cite{McCann2019}.

We provide the details in Section \ref{s:supermartingale}.

\subsubsection{Harmonic transport}

Another example of an application of our theory concerns harmonic transports. Suppose that $\mathcal{A}$ is the space of harmonic functions on an open set $\Omega\subset\mathbb{R}^n$. 
By Theorem \ref{thm:holoin} there exists a non-negative, continuous, proper function $p$ in the complete lattice cone $\mathcal{H}$ generated by $\mathcal{A}$, so that Theorem \ref{thm:strassenINTRO} is applicable. We shall consider two measures $\mu,\nu$, integrable with respect to $p$, and such that $\mu\prec_{\mathcal{H}}\nu$. 

If we assumed that $\mu$ and $\nu$ are in order with respect to the cone $\mathcal{S}$ of all upper semi-continuous subharmonic functions, then the resolution of the Skorokhod embedding problem, see e.g.  \cite{Palmer2019}, would show that there exists a stopping time $\tau$ and a Brownian motion $(B_t)_{t\in [0,\infty)}$ such that $B_0$ and $B_{\tau}$ are distributed according to $\mu$ and $\nu$ respectively.



In  Section \ref{s:harmonic} we discuss in detail how Theorem \ref{thm:strassenINTRO} and Theorem \ref{thm:partitionINTRO} can be applied to provide, given $\mu\prec_{\mathcal{H}}\nu$, a partitioning of $\Omega$ into $\mathcal{A}$-convex sets, that yields constraints for all $\mathcal{H}$-transports between $\mu$ and $\nu$.
We show an example, Example \ref{exa:harmonic}, of a pair of measures for which such partitioning is strictly finer than the associated convex partitioning.








\subsubsection{Relaxed martingale transport}

In \cite{Guo2019}, Guo and Ob\l\'oj introduced a relaxed notion of approximate martingale transport in order to develop a numerical scheme for one-dimensional, multi-step optimal martingale transport problem.
 An important issue related to this is the continuity of the martingale optimal transport, studied e.g. by Wiesel in  \cite{Wiesel2023}.


In the setting of one-step martingales, given $\delta\geq 0$ and two probability measures $\mu,\nu$ on the real  line in convex order with finite first moments, one looks for transports $\pi$ between $\mu$ and $\nu$ with the gap from a genuine martingale transport between $\mu$ and $\nu$ controlled by $\delta$.
More specifically, we say that a transport $\pi$ between $\mu$ and $\nu$ is a $\delta$-approximating martingale transport provided that for all affine $a\in\mathcal{A}$ and all non-negative, bounded and measurable $h$
\begin{equation*}
\int_{\mathbb{R}\times\mathbb{R}}\big(h(\omega_1)a(\omega_1)-h(\omega_1)a(\omega_2)\big)\, d\pi(\omega_1,\omega_2)\leq \delta\norm{h}_{\mathcal{D}_1(\mathbb{R})}\norm{a}_{\mathcal{D}_{p+1}(\mathbb{R})},
\end{equation*}
where $p$ denotes a norm on the real line. 
Clearly, if $\mu\prec_{\mathcal{F}}\nu$ then the set of such transports is non-empty, by the Strassen theorem, cf. Theorem \ref{thm:strassenINTRO}.

More generally, if we have a lattice cone $\mathcal{F}$ generated by a set $\mathcal{G}$ of functions on $\Omega$, we may  look at the set $\Gamma_{\mathcal{F},\delta}(\mu,\nu)$ 
 of all transports between $\mu$ and $\nu$ such that for all $g\in\mathcal{G}$ and all bounded, non-negative measurable $h$
\begin{equation}\label{eqn:ineqa}
\int_{\Omega\times\Omega}\big(h(\omega_1)g(\omega_1)-h(\omega_1)g(\omega_2)\big)\, d\pi(\omega_1,\omega_2)\leq \delta\norm{h}_{\mathcal{D}_1(\Omega)}\norm{g}_{\mathcal{D}_{p+1}(\Omega)}.
\end{equation}
Again, if $\mu\prec_{\mathcal{F}}\nu$, then such inequality is satisfied by any transport $\Gamma_{\mathcal{F}}(\mu,\nu)$, which is non-empty by the Strassen theorem.

We claim that for any  $\delta> 0$, the condition (\ref{eqn:ineqa}) does not impose any constraints on the supports of disintegrations of a $\delta$-approximating martingale transport, unlike in the case $\delta=0$. This provides a negative resolution of a question of Cox, \cite{Cox2022}. 

We provide a characterisation of polar sets with respect to $\delta$-approximating $\mathcal{F}$-transports, when $\delta>0$. We demonstrate that these sets coincide with the sets that are polar with respect to all transports between $\mu$ and $\nu$.

We provide the precise definition and the statement of the theorem in Section \ref{s:relaxed}, Definition \ref{def:delta} and Theorem \ref{thm:relaxed}.

\subsubsection{Finite-dimensional approximations and infinite-dimensional setting}

In Theorem \ref{thm:partitionINTRO} we assume that all $\mathcal{F}$-transports between $\mu$ and $\nu$ are local. Therefore, we do not immediately obtain a partitioning if there exists a non-local $\mathcal{F}$-transport between given two measures in $\mathcal{F}$-order.

A resolution of that obstacle can be achieved by the procedure of finite-dimensional approximations.

\begin{theorem}\label{thm:apircintro}
    Suppose that $\mu\prec_{\mathcal{F}}\nu$. Suppose that $\mathcal{G}_0$ be at most countable subset of $\mathcal{D}_{p+1}(\Omega)$. Then there exist an assignment to any $\omega\in\Omega$ of an $\mathcal{A}$-convex set $\mathrm{apirc}_{\mathcal{G}_0}(\mu,\nu)(\omega)\subset\mathcal{A}^*$  and  a Borel measurable set $B\subset\Omega$ with $\mu(B)=1$ such that whenever $\omega_1,\omega_2\in B$  then
    \begin{equation*}
        \mathrm{apirc}_{\mathcal{G}_0}(\mu,\nu)(\omega_1)\cap \mathrm{apirc}_{\mathcal{G}_0}(\mu,\nu)(\omega_2)=\emptyset
    \end{equation*}
    or
    \begin{equation*}
        \mathrm{apirc}_{\mathcal{G}_0}(\mu,\nu)(\omega_1)= \mathrm{apirc}_{\mathcal{G}_0}(\mu,\nu)(\omega_2).
    \end{equation*}
    Moreover, the sets $\mathrm{apirc}_{\mathcal{G}_0}(\mu,\nu)(\omega)$ are convex and for any  $\lambda\in\Lambda_{\mathcal{F}}(\mu,\nu)$
    \begin{equation*}
        \mathrm{supp}\lambda(\omega,\cdot)\subset \Phi^{-1}\big(\mathrm{cl}\big(\mathrm{apirc}_{\mathcal{G}_0}(\mu,\nu)(\omega)\big)\big)\text{  for }\mu\text{-almost every }\omega\in\Omega.
    \end{equation*}
    Furthermore, for $\mathcal{B}=\mathrm{Cone}\mathcal{G}_0\cap(-\mathrm{Cone}\mathcal{G}_0)$,
    \begin{equation*}
       \Phi_{\mathcal{B}}(\omega)\in R_{\mathcal{B}}\big(\mathrm{cl}\big(\mathrm{apirc}_{\mathcal{G}_0}(\mu,\nu)(\omega)\big)\big)\text{ for }\mu\text{-almost every }\omega\in\Omega.
        \end{equation*}
    In particular, if $\mathcal{G}_0$ is symmetric, then 
    \begin{equation*}
        \Phi(\omega)\in \mathrm{cl}\big(\mathrm{apirc}_{\mathcal{G}_0}(\mu,\nu)(\omega)\big)\text{ for }\mu\text{-almost every }\omega\in\Omega.
    \end{equation*}
\end{theorem}

The sets $\mathrm{apirc}_{\mathcal{G}_0}(\mu,\nu)(\cdot)$ shall be called the $\mathcal{G}_0$-\emph{approximate irreducible components}. As we shall see in  Definition \ref{def:approximate}, they depend only on $\mu,\nu,\mathcal{G}_0,\mathcal{G}$.

The difference in conclusions of Theorem \ref{thm:apircintro} and Theorem \ref{thm:partitionINTRO} is that in the  former theorem we do not claim that the obtained partitioning is the finest, unlike in Theorem \ref{thm:partitionINTRO}. The fact that the irreducible components of Theorem \ref{thm:partitionINTRO} are the smallest convex sets that satisfy the stipulated conditions allows for a characterisation of the polar sets in  Theorem \ref{thm:polarinto}.

Let us briefly describe the idea of the proof. For any finite family $Z\subset\mathcal{G}_0$, we let $\mathcal{A}_Z$ to be the linear span of $Z$ and $\mathcal{F}_Z$ to be the lattice cone generated by $\mathcal{G}$. Since $\mathcal{F}_Z\subset\mathcal{F}$, we see that $\mu\prec_{\mathcal{F}}\nu$ implies that $\mu\prec_{\mathcal{F}_Z}\nu$. Let $\epsilon>0$.
  Let $\Gamma_{\mathcal{F}_Z}(\mu,\nu)$ be the non-empty set of all $\mathcal{F}_Z$-transports between $\mu$ and $\nu$ and let $\Lambda_{\mathcal{F}_Z}(\mu,\nu)$ be the corresponding set of all disintegrations of elements of $\Gamma_{\mathcal{F}_Z}(\mu,\nu)$ with respect to the first co-ordinate, $\Phi_{\mathcal{A}_Z}$ be the assignment to $\omega\in\Omega$ the corresponding evaluation functional on $\mathcal{A}_Z$.
 The closed convex hulls of 
 \begin{equation*}
     \Phi_{\mathcal{A}_Z}(\mathrm{supp}\lambda(\omega,\cdot))\subset\mathcal{A}_Z^*
 \end{equation*} 
 are finite-dimensional, where $\lambda\in\Lambda_{\mathcal{F}}(\mu,\nu)$ is a maximal disintegration of $\Gamma_{\mathcal{F}}(\mu,\nu)$. We may employ Theorem \ref{thm:partitionINTRO}, to infer that the sets $\mathrm{eirc}_{\mathcal{A}_Z}(\mu,\nu)(\omega)$ of all extensions of the functionals in $\mathrm{irc}_{\mathcal{A}_Z}(\mu,\nu)(\omega)$ to $\mathcal{A}$, for $\omega\in\Omega$, form a partitioning of $\Omega$, up to a set $B_Z\subset\Omega$ of $\mu$-measure zero.
Then the sets
\begin{equation*}
    \bigcap\Big\{\mathrm{eirc}_{\mathcal{A}_Z}(\mu,\nu)(\omega)\mid Z\subset \mathcal{G}_0\text{ is finite}\Big\}\text{, for }\omega\in\Omega,
\end{equation*}
form a partitioning of $\Omega$ up to a set of $\mu$-measure zero.


We provide further details in Section \ref{s:approx}. An application of the theorem is presented in Section \ref{s:optimal}, where we consider an application of the general localisation scheme to the optimal transport localisation. In Section \ref{s:harmonic} we apply the finite-dimensional approximations approach to the setting of harmonic transports, in Section \ref{s:martingaleproblem} we use the approximations to provide a result for the solutions of the martingale problems.

\subsubsection{Martingale problems}

Yet another example pertains to the solutions of martingale problems. Suppose that $L$ is a linear operator on the subspace $\mathcal{D}$ of the space of Borel functions, $\mu$ is a Borel probability measure and $(X_t)_{t\in [0,1]}$ is a solution of the martingale problem for $(L,\mu)$, i.e., for any $f\in\mathcal{D}$
\begin{equation*}
    f(X_t)-\int_0^t(Lf)(X_s)\,ds
\end{equation*}
is a well-defined martingale. Let $\mathcal{A}$ denote the kernel of $L$. Let $\mathcal{G}_0\subset\mathcal{A}$ be a countable subset and let $\nu$ be the distribution of $X_1$. Then there is a partitioning of the underlying space into $\mathcal{A}$-convex sets $\Phi^{-1}(\mathrm{apirc}_{\mathcal{G}_0}(\mu,\nu)(\cdot))$ such that for all $t\in [0,1]$
\begin{equation*}
    X_t\in \Phi^{-1}\big(\mathrm{cl}(\mathrm{apirc}_{\mathcal{G}_0}(\mu,\nu)(X_0))\big)\text{ almost surely.}
\end{equation*}
As illustrated in Section \ref{s:harmonic}, the obtained partitioning can be non-trivial.

These ideas are discussed in a greater detail in Section \ref{s:martingaleproblem}.

Let us mention that there are no restrictions as to what the subspace $\mathcal{A}$ can be. We may choose $\mathcal{A}$ to consist of real parts of holomorphic functions. This choice would lead to holomorphically convex paving for any holomorphic martingale.

\section{Preliminaries}\label{s:prelim}

Let us recall several definitions.



\begin{definition}
Let $\Omega$, $Z$ be topological spaces. A map $p\colon \Omega\to Z$ is said to be proper whenever preimages $p^{-1}(Z)$ of compact sets in $Z$ are compact sets in $\Omega$.
\end{definition}

\begin{definition}
Let $\Omega$ be a set and let $\mathcal{A}$ be a set of functions on $\Omega$. The coarsest topology on $\Omega$ with respect to which all functions in $\mathcal{A}$ are continuous we shall call the topology generated by $\mathcal{A}$ and denote by $\tau(\mathcal{A})$.
\end{definition}

\begin{definition}
Let $\Omega$ be a set and let $p$ be a function on $\Omega$ bounded from below by a positive constant. We shall say that a function $a$ on $\Omega$ is of $p$-growth if there exists $C>0$ such that for all $\omega\in\Omega$ there is
\begin{equation*}
\abs{a(\omega)}\leq Cp(\omega).
\end{equation*}
We denote by $\mathcal{D}_p(\Omega)$ the space of all functions of $p$-growth, with norm 
\begin{equation*}
\norm{a}_{\mathcal{D}_p(\Omega)}=\sup\bigg\{\frac{\abs{a(\omega)}}{p(\omega)}\mid \omega\in\Omega\bigg\}.
\end{equation*}
If $\Omega$ is a topological space, then we let $\mathcal{C}_p(\Omega)$ denote the space of all continuous functions of $p$-growth.
We let $\mathcal{C}(\Omega)$ denote the space of all continuous functions on $\Omega$ and $\mathcal{M}_p(\Omega)$ denote the space of signed Borel measures $\mu$ on $\Omega$ with 
\begin{equation*}
\int_{\Omega}p\, d\abs{\mu}<\infty.
\end{equation*}
For $\mu\in  \mathcal{M}_p(\Omega)$ we  set 
\begin{equation*}
\norm{\mu}_{\mathcal{M}_p(\Omega)}=\int_{\Omega}p\, d\abs{\mu}.
\end{equation*}
We let $\mathcal{P}_p(\Omega)$ denote the subset of $\mathcal{M}_p(\Omega)$ consisting of probability measures.
\end{definition}

\begin{definition}\label{def:baire}
Let $\Omega$ be a topological space, with topology $\tau$. Then:
\begin{enumerate}
\item $\sigma(\tau)$ is called the Borel $\sigma$-algebra on $\Omega$,
\item a measure $\mu$ on $\sigma(\tau)$ is called a Radon measure if for any set $B\in\sigma(\tau)$ and $\epsilon>0$ there exists a compact set $K\subset B$ such that $\abs{\mu}(B\setminus K)<\epsilon$,
\end{enumerate}
\end{definition}

\begin{definition}\label{def:prelim}
Let $\Omega$ be a set. Let $\mathcal{G}$ be a set of functions on $\Omega$ with values in $(-\infty,\infty]$. A set $\mathcal{K}$ of functions on $\Omega$ is said to be:
\begin{enumerate}
\item a convex cone whenever $\alpha f+\beta g \in\mathcal{K}$ for any $f,g\in\mathcal{K}$ and any numbers $\alpha,\beta\geq 0$;
\item stable under maxima provided that the maximum $f\vee g$ of any two functions   $f,g\in\mathcal{K}$ belongs to $\mathcal{K}$;
\item stable under suprema provided that the supremum $\sup\{f_i\mid i \in I\}$  of any family $(f_i)_{i\in I}\subset\mathcal{K}$ belongs to $\mathcal{K}$;
\item a lattice cone whenever it is a convex cone stable under maxima;
\item the lattice cone generated by $\mathcal{G}$ whenever it is the smallest lattice cone of functions on $\Omega$ containing $\mathcal{G}$;
\item a complete lattice cone whenever it is a convex cone stable under suprema;
\item the complete lattice cone generated by $\mathcal{G}$ whenever it is the smallest complete lattice cone of functions on $\Omega$ containing $\mathcal{G}$;
 \item separating points of $\Omega$ whenever for any two distinct $\omega_1,\omega_2\in\Omega$, there exist $f,g\in\mathcal{K}$ such that $f(\omega_1)\neq g(\omega_2)$.
\end{enumerate}
\end{definition}

Let us stress that we allow the functions to take value $+\infty$. Note however that this is not allowed for linear subspaces.

Let us also recall several definitions from \cite{Ciosmak20242}. For a complete description we refer the reader therein.

\begin{definition}
    Suppose that $\mu,\nu\in\mathcal{P}_p(\Omega)$ are two Radon probability measures.
Any Radon probability measure $\pi$ on $\Omega\times\Omega$ with marginals $\mu,\nu$ and such that for all $g\in\mathcal{G}$ and all bounded, non-negative measurable $h$ we have
\begin{equation*}
\int_{\Omega}g(\omega_1)h(\omega_1)\,d\mu(\omega_1)\leq \int_{\Omega\times\Omega}g(\omega_2)h(\omega_1)\, d\pi(\omega_1,\omega_2)
\end{equation*}
we shall call an $\mathcal{F}$-transport between $\mu$ and $\nu$. 
We shall denote the set of all such measures by $\Gamma_{\mathcal{F}}(\mu,\nu)$. 
\end{definition}

\begin{definition}
Let $\pi\in\Gamma_{\mathcal{F}}(\mu,\nu)$. Let $T(\pi)$ be the set of all disintegrations
\begin{equation*}
\lambda(\cdot,\cdot)\colon  \Omega\times \sigma(\tau(\mathcal{A}))\to\mathbb{R}
\end{equation*}
of $\pi$ with respect to the projection $p_1$ on the first co-ordinate. That is,
\begin{enumerate}
\item for every $\omega\in \Omega$, $ \sigma(\tau(\mathcal{A}))\ni B\mapsto \lambda(\omega,B)\in\mathbb{R}$ is a Radon probability measure on $\Omega$,
\item for every $B\in\sigma(\tau(\mathcal{A}))$ the map $\Omega\ni \omega\mapsto\lambda(\omega,B)\in\mathbb{R}$ is Borel measurable,
\item for all Borel sets $B\subset\Omega\times\Omega$ and $E\subset\Omega$, 
\begin{equation*}
\pi\big(B\cap r^{-1}(E)\big)=\int_E\lambda(\omega,B_{\omega})\, d\mu(\omega),
\end{equation*}
where 
\begin{equation*}
B_{\omega}=\{\omega'\in\Omega\mid (\omega,\omega')\in\Omega\}.
\end{equation*}
\end{enumerate}
\end{definition}

\begin{definition}
We shall denote by $\Lambda_{\mathcal{F}}(\mu,\nu)$ the set of all maps $\lambda\in T(\pi)$ for some $\pi\in \Gamma_{\mathcal{F}}(\mu,\nu)$.
\end{definition}

\begin{definition}
    A Borel set $U\subset\Omega\times\Omega$ is said to be a polar set with respect to $\Gamma_{\mathcal{F}}(\mu,\nu)$ if $\pi(U)=0$ for all $\pi\in\Gamma_{\mathcal{F}}(\mu,\nu)$. Similarly, if $\Gamma\subset \Gamma(\mu,\nu)$, then $U$ is said to be a polar with respect to $\Gamma$ if $\pi(U)=0$ for all $\pi\in\Gamma$.
\end{definition}

\begin{definition}\label{def:maxdis}
A map $\lambda\in \Lambda_{\mathcal{F}}(\mu,\nu)$ such that for any $\lambda'\in \Lambda_{\mathcal{F}}(\mu,\nu)$ there is
\begin{equation*}
\mathrm{supp}\lambda'(\omega,\cdot)\subset\mathrm{supp}\lambda(\omega,\cdot)\text{ for }\mu\text{-almost every }\omega\in \Omega,
\end{equation*}
we shall call a maximal disintegration of $\Gamma_{\mathcal{F}}(\mu,\nu)$ with respect to $p_1$, or, briefly, a maximal disintegration of  $\Gamma_{\mathcal{F}}(\mu,\nu)$. Note here that we consider the supports of measures with respect to the topology $\tau(\mathcal{G})$.
\end{definition}

\begin{definition}
    Let $\mathcal{A}\subset\mathcal{D}_p(\Omega)$. We define the Gelfand transform $\Phi\colon\Omega\to\mathcal{A}^*$ by the formula
    \begin{equation*}
        \Phi(\omega)(a)=a(\omega)\text{ for all }\omega\in\Omega\text{ and }a\in\mathcal{A}.
    \end{equation*}
\end{definition}

\begin{definition}\label{def:local}
    A transport $\pi$ between $\mu$ and $\nu$ is called local whenever for $\mu$-almost every $\omega\in\Omega$ the set
\begin{equation*}
    \mathrm{clConv}_H\Phi(\mathrm{supp}\lambda(\omega,\cdot))\subset\mathcal{A}^*
\end{equation*}
is a finite-dimensional convex set. Here $\lambda\in T(\pi)$ is a disintegration of $\pi$ with respect to the first co-ordinate, and $H$ denotes the complete lattice cone of convex, lower semi-continuous functions on $\mathcal{A}^*$, cf. Section \ref{s:fconvex}.
\end{definition}

\begin{definition}
    If $\mu\prec_{\mathcal{F}}\nu$, then for $\omega\in\Omega$ we define
    \begin{equation*}
        \mathrm{irc}_{\mathcal{A}}(\mu,\nu)(\omega)=\mathrm{rint}\Phi(\mathrm{clConv}_H\mathrm{supp}\lambda(\omega,\cdot)),
    \end{equation*}
    where $\lambda\in\Lambda_{\mathcal{F}}(\mu,\nu)$ is a maximal disintegration of $\Gamma_{\mathcal{F}}(\mu,\nu)$. The sets $\mathrm{irc}_{\mathcal{A}}(\mu,\nu)(\cdot)$ we shall call irreducible components. 
\end{definition}

\section{$\mathcal{F}$-convex sets}\label{s:fconvex}

Throughout this section we shall consider a set $\Omega$ and a linear space $\mathcal{A}$ of functions on $\Omega$. We shall denote by $\mathcal{F}$ the complete lattice cone of functions on $\Omega$ generated by $\mathcal{A}$. 

In \cite[Section 4]{Ciosmak20242} we have proven a general version of the Strassen theorem. We assumed that any function in $\mathcal{A}$ is of $p$-growth for some proper, non-negative function $p$ in the complete lattice cone $\mathcal{H}$ generated by $\mathcal{A}$.

We shall now cite a theorem of \cite{Ciosmak2023} that gives necessary and sufficient conditions for existence of such a function.

If $\Omega$ is a set and $\mathcal{A}$ is a linear space of functions on $\Omega$, we shall say that a net $(\omega_{\alpha})_{\alpha\in A}$ is a Cauchy net in $\Omega$ whenever for each $a\in\mathcal{A}$, $(a(\omega_{\alpha}))_{\alpha\in A}$ is a Cauchy net in $\mathbb{R}$. We shall say that $\Omega$ is complete with respect to $\mathcal{A}$ whenever every Cauchy net in $\Omega$ is convergent in $\tau(\mathcal{A})$. We refer to \cite{Schaefer1971} for a short background on Cauchy nets and completeness in relation to linear topological spaces.

\begin{theorem}\label{thm:holoin}
Let $\Omega$ be a $\sigma$-compact and locally compact Hausdorff topological space.  Let $\mathcal{A}$ be a linear space of continuous functions on $\Omega$ that contains constants and separates points of $\Omega$.
Then the following conditions are equivalent:
\begin{enumerate}
\item\label{i:fconstantin} for any compact set $K\subset\Omega$ and for any $C\geq 1$ the set
\begin{equation*}
\{\omega\in\Omega\mid a(\omega)\leq C\sup \abs{a}(K)\text{ for all }a\in\mathcal{A}\}
\end{equation*}
is a compact subset of $\Omega$,
\item\label{i:exhaustion} there exists a non-negative, continuous and proper function $p$ on  $\Omega$ such that  
\begin{equation*}
p=\sup\{a_{\alpha}\mid \alpha\in A\}
\end{equation*}
for some symmetric family $(a_{\alpha})_{\alpha\in A}$
of functions in $\mathcal{A}$,
\item\label{i:exhaustionbound} there exists a non-negative, proper function $p$ on  $\Omega$, bounded on compactae, such that  
\begin{equation*}
p=\sup\{a_{\alpha}\mid \alpha\in A\}
\end{equation*} for some symmetric family $(a_{\alpha})_{\alpha\in A}$ of functions in $\mathcal{A}$,
\item\label{i:compactae}
there exists a family $(P_i)_{i=1}^{\infty}$ of compact, symmetric $\mathcal{A}$-polygons such that
\begin{equation*}
\bigcup_{i=1}^{\infty}P_i=\Omega, P_i\subset\mathrm{int}P_{i+1}\text{  for }i=1,2,\dotsc.
\end{equation*} 
\end{enumerate}
Moreover, the suprema in  \ref{i:exhaustion} and in \ref{i:exhaustionbound} may be locally taken over finite families of functions in $\mathcal{A}$.

Suppose moreover that $\Omega$ is separable and $\mathcal{A}$ is closed with respect to the compact-open topology. Then the above conditions are equivalent to $\Omega$ being complete with respect to $\mathcal{A}$.
If moreover $\mathcal{A}$ consists of real-parts of a complex algebra, then the above conditions are equivalent to $\Omega$ being an $\mathcal{A}$-space.
\end{theorem}


Let us also note that another theorem of \cite{Ciosmak2023} is concerned with charactersation of the existence of exhaustion functions belonging to a cone of functions satisfying the maximum principle.

\begin{definition}
Let $\mathcal{F}$ be a cone of functions on topological space $\Omega$. Let $S\subset\Omega$. We define closed $\mathcal{F}$-convex hull of $S$ by the formula
\begin{equation*}
\mathrm{clConv}_{\mathcal{F}}S=\{\omega\in\Omega\mid f(\omega)\leq \sup f(S)\text{ for all }f\in\mathcal{F}\}.
\end{equation*}
\end{definition}

The lemma and the proposition below are proven in \cite{Ciosmak2023}.

\begin{lemma}\label{lem:convex}
Suppose that $\mathcal{F}$ is a complete lattice cone generated by a linear subspace $\mathcal{A}$. Let $S\subset\Omega$. Then
\begin{equation*}
\mathrm{clConv}_{\mathcal{F}}S=\{\omega\in\Omega\mid a(\omega)\leq \sup a(S)\text{ for all }a\in\mathcal{A}\}.
\end{equation*}
\end{lemma}

\begin{proposition}\label{pro:convex}
Let $S\subset\Omega$. Then
\begin{equation*}
\mathrm{clConv}_{\mathcal{F}}S=\Phi^{-1}(\mathrm{clConv}_G\Phi(S)),
\end{equation*}
where $G$ is the cone of convex functions on $\mathcal{A}^*$, lower semi-continuous with respect to $\sigma(\tau(\mathcal{A}))$-topology.
\end{proposition}

The above lemma and proposition allowed in \cite{Ciosmak2023} for a general definition of $\mathcal{F}$-convex set.

\begin{definition}\label{def:conv}
A set $S\subset\Omega$ is said to be $\mathcal{F}$-convex whenever $S= \Phi^{-1}(\mathrm{Conv}_G\Phi(S))$.
\end{definition}

\section{Martingale transports in infinite dimensions}\label{s:martingale}

The example that directly generalises the results of \cite{Juillet2016, Nutz2017, Touzi2019, Obloj2017} concerns the martingale transport of Radon probability measures on infinite-dimensional spaces. 

Let $X^*$ be an infinite dimensional normed space $X^*$, that is dual to some normed space $X$. We equip it with the weak* topology. The set of all convex,  lower semi-continuous with  respect to the weak* topology on $X^*$, functions $f\colon X^*\to\mathbb{R}$ shall be denoted by $\mathcal{H}$. It is the complete lattice cone generated by the linear space $\mathcal{A}$ of all weakly* continuous affine functions on $X^*$. Let $\mathcal{F}$ be the lattice cone generated by $\mathcal{A}$. 


\begin{theorem}\label{thm:martingale}
Suppose that $X^*$ is separable Banach space, dual to $X$. Let $\mu,\nu\in\mathcal{P}(\Omega)$ be two Radon measures such that 
\begin{equation*}
    \int_{X^*}\norm{\cdot}\, d\mu<\infty,   \int_{X}\norm{\cdot}\, d\nu<\infty
\end{equation*}
and $\mu\prec_{\mathcal{F}}\nu$.
Suppose that any $\mathcal{F}$-transport between $\mu$ and $\nu$ is local.
Then there exists a Borel measurable set $B\subset X^*$ with $\mu(B)=1$ such that whenever $x_1^*,x_2^*\in B$ then 
\begin{equation*}
\mathrm{irc}_{\mathcal{A}}(\mu,\nu)(x_1^*)\cap\mathrm{irc}_{\mathcal{A}}(\mu,\nu)(x_2^*)=\emptyset
\end{equation*}
or
\begin{equation*}
\mathrm{irc}_{\mathcal{A}}(\mu,\nu)(x_1^*)=\mathrm{irc}_{\mathcal{A}}(\mu,\nu)(x_2^*).
\end{equation*}
These sets are relatively open, convex and finite-dimensional subsets of $X^*$.
Moreover, for any disintegration $\lambda\in\Lambda_{\mathcal{F}}(\mu,\nu)$
\begin{equation*}
    \mathrm{supp}\lambda(\omega,\cdot)\subset \mathrm{cl}\big(\mathrm{irc}_{\mathcal{A}}(\mu,\nu)(x^*)\big)\text{ for }\mu\text{-almost every }x^*\in X^*.
\end{equation*}
The sets $\mathrm{irc}_{\mathcal{A}}(\mu,\nu)(\cdot)$ are the smallest convex sets that satisfy the above condition.
Furthermore, 
\begin{equation*}
    \omega\in \mathrm{irc}_{\mathcal{A}}(\mu,\nu)(x^*)\text{ for }\mu\text{-almost every }x^*\in X^*.
\end{equation*}
    
Let $U\subset\Omega\times\Omega$ be a Borel set such that for $\mu $-almost every $x^*\in X^*$
\begin{equation*}
U_{x^*}\subset\mathrm{irc}_{\mathcal{A}}(\mu,\nu)(x^*).
\end{equation*}
Then $U$ is $\Gamma_{\mathcal{F}}(\mu,\nu)$-polar set if and only if there exist Borel sets $N_1,N_2\subset\Omega$ with
\begin{equation*}
\mu(N_1)=0,\nu(N_2)=0
\end{equation*}
and 
\begin{equation*}
U\subset (N_1\times\Omega)\cup (\Omega\times N_2).
\end{equation*} 
\end{theorem}
\begin{proof}
    The theorem is a direct application of Theorem \ref{thm:partitionINTRO} and Theorem \ref{thm:polarinto}. We need to verify the assumptions of these theorems. 
Let us observe that by the Banach--Alaoglu theorem the norm of $X^*$ is proper. It belongs to $\mathcal{H}$ thanks to the Hahn--Banach theorem. Clearly, $\mathcal{A}$ separates points of $X^*$ and consists of functions of $(p+1)$-growth.
Therefore Theorem \ref{thm:strassenINTRO} shows that the set $\Gamma_{\mathcal{F}}(\mu,\nu)$ of martingale transports between $\mu$ and $\nu$ is non-empty.

The assumption tells us that any $\mathcal{F}$-transport between $\mu$ and $\nu$ is local.
Theorem \ref{thm:partitionINTRO} shows that there exists the finest partitioning of $\mathcal{A}^*$, up to a set of  $\mu$-measure zero, into relatively open convex and finite-dimensional sets $\mathrm{irc}_{\mathcal{F}}(\mu,\nu)(\cdot)$. Since any element is assumed to give value $1$ to a constant function $1$, and $\mathcal{A}$ can be identified with $X+\mathbb{R}1$, we see that the sets can be identified with subsets of $X^*$. 

Still by Theorem \ref{thm:partitionINTRO}, any for $\mathcal{F}$-transport and any its disintegration $\lambda\in\Lambda_{\mathcal{F}}(\mu,\nu)$, 
\begin{equation*}
    \mathrm{supp}\lambda(x^*,\cdot)\subset \mathrm{cl}\,\mathrm{irc}_{\mathcal{F}}(\mu,\nu)(x^*)\text{ for }\mu\text{-almost every  }x^*\in X^*,
\end{equation*}
and
\begin{equation*}
 x^*\in\mathrm{irc}_{\mathcal{F}}(\mu,\nu)(x^*)\text{ for }\mu\text{-almost every  }x^*\in X^*.
\end{equation*}
The characterisation of polar sets follows by Theorem \ref{thm:polarinto}.
\end{proof}

The characterisation of polar sets gives an affirmative answer to a generalisation of a conjecture of Ob\l\'oj and Siorpaes \cite{Obloj2017} to the infinite-dimensional setting. It says, in a short formulation, that a set $U$ that has sections contained in the corresponding irreducible components, is a polar set if and only if it is polar with respect to all transports between $\mu$ and $\nu$.

Measure $\mu$ may be disintegrated, thanks to \cite[Theorem 13.16]{Ciosmak20242}, in such a way that any conditional measure admits a single irreducible component.

If the considered transports are not assumed to be local, then we may still provide a partitioning, thanks to Section \ref{s:approx}, Theorem \ref{thm:apir}. We do not however obtain a characterisation of polar sets.

The following gives also a characterisation of intersections of the closures of the irreducible components.
\begin{theorem}\label{thm:intersections}
   Suppose that $X^*$ is separable Banach space, dual to $X$. Let $\mu,\nu\in\mathcal{P}(\Omega)$ be two Radon measures such that 
\begin{equation*}
    \int_{X^*}\norm{\cdot}\, d\mu<\infty,   \int_{X}\norm{\cdot}\, d\nu<\infty
\end{equation*}
and $\mu\prec_{\mathcal{F}}\nu$.
Suppose that any $\mathcal{F}$-transport between $\mu$ and $\nu$ is local.
Then there exists a Borel measurable set $B\subset X^*$ with $\mu(B)=1$ such that whenever $x_1^*,x_2^*\in B$, and if 
\begin{equation*}
    x^*\in \mathrm{cl}\big(\mathrm{irc}_{\mathcal{A}}(\mu,\nu)(x_1^*)\big)\cap\mathrm{cl}\big(\mathrm{irc}_{\mathcal{A}}(\mu,\nu)(x_2^*)\big)
\end{equation*}
is  such that for $i=1,2$ there exists $\eta_i\in\mathcal{P}_{\norm{\cdot}}(\Omega)$ and $C_i>0$ such that $\eta_i\leq C_i\lambda(\omega_i,\cdot)$, where $\lambda\in\Lambda_{\mathcal{F}}(\mu,\nu)$ is a maximal disintegration, and
\begin{equation*}
x^*(x)=\int_{\Omega}x\, d\eta_i\text{ for all }x\in\mathcal{X},
\end{equation*}
then we have the equality of the Gleason parts
\begin{equation*}
G\big(x^*,\mathrm{cl}\big(\mathrm{irc}_{\mathcal{A}}(\mu,\nu)(x_1^*)\big)\big)=G\big(x^*,\mathrm{cl}\big(\mathrm{irc}_{\mathcal{A}}(\mu,\nu)(x_2^*)\big)\big).
\end{equation*}
\end{theorem}
\begin{proof}
    The theorem is a consequence of Theorem \ref{thm:finestrucfinintro}.
\end{proof}

\section{Submartingale transports in finite and infinite dimensions}\label{s:supermartingale}

Let us now consider the example that generalises the results of \cite{Nutz2018} to any  submartingale transports in finite dimensions and to local submartingale transports in infinite dimensions. Let us note that \cite{Nutz2018} is concerned with supermartingale transports -- these however correspond to submartingale transports under flipping of the signs in the underlying space.

Let $X^*$ be an infinite dimensional normed space $X^*$, that is dual to some normed space $X$. We equip it with the weak* topology. 
Let us pick a subset $\mathcal{G}_0$ of $X$, whose linear span is dense in $X$. Each element of $\mathcal{G}_0$ defines a continuous linear function on $X^*$. We take $\mathcal{G}$ to be the union of $\mathcal{G}_0$ and constant functions. Let $\mathcal{F}$ be the lattice cone generated by $\mathcal{G}$.

Then, see \cite[Lemma 2.6]{Ciosmak20242}, the complete lattice cone $\mathcal{H}$ generated by $\mathcal{G}$ is the set of lower semi-continuous, convex functions that are non-decreasing in the directions of the polar cone to the cone spanned by $\mathcal{G}$, i.e., 
\begin{equation*}
\mathcal{G}^*=\{x^*\in X^*\mid x^*(g)\geq 0\text{ for all }g\in\mathcal{G}\}.
\end{equation*}
This is to say, for $h\in\mathcal{H}$, 
\begin{equation*}
  h(x^*)\leq   h(x^*+g^*)\text{ for any }x^*\in X^*\text{ and }g^*\in\mathcal{G}^*.
\end{equation*}

\begin{remark}\label{rem:submart}
Let us note that when $\mathcal{G}_0$ is symmetric, i.e., $\mathcal{G}_0=-\mathcal{G}_0$, then the set $\mathcal{H}$ is the cone of lower semi-continuous, convex functions on $X^*$. In this case, all the considered transports will be martingale transports. Consider the conical hull $\mathrm{Cone}\mathcal{G}$ of $\mathcal{G}$. It can be written as 
\begin{equation*}
    \mathrm{Cone}\mathcal{G}=\mathcal{G}_1+ \mathcal{G}_2,
\end{equation*}
where $\mathcal{G}_1$ is the maximal linear subspace $\mathcal{G}_1\subset \mathrm{Cone}\mathcal{G}$ and $\mathcal{G}_2$ is a convex cone. Clearly, 
\begin{equation*}
    \mathcal{G}_1=\mathrm{Cone}\mathcal{G}\cap (-\mathrm{Cone}\mathcal{G}).
\end{equation*} 
Then any $\mathcal{F}$-transport is a martingale transport in directions of $\mathcal{G}_1$, and submartingale transport in directions of $\mathcal{G}_2$. That is, if $\pi\in\Gamma_{\mathcal{F}}(\mu,\nu)$, then for all $g\in\mathcal{G}_1$ we have
\begin{equation*}
    \int_{\Omega}h(\omega_1)g(\omega_1)=\int_{\Omega\times\Omega}h(\omega_1)g(\omega_2)\, d\pi(\omega_1,\omega_2)
\end{equation*}
for all bounded, non-negative and measurable $h$. If $g\in\mathcal{G}_2$, we know only that the left hand-side of the above equality is bounded by its right hand-side.
\end{remark}

By Remark \ref{rem:submart}, the following theorem covers the case of martingale-submartingale transports, i.e., transports that are martingale measures in certain directions and that are submartingale measures in certain directions.

\begin{theorem}\label{thm:submartingale}
Suppose that $X^*$ is separable Banach space, dual to $X$. Let $\mu,\nu\in\mathcal{P}(\Omega)$ be two Radon measures such that 
\begin{equation*}
    \int_{X^*}\norm{\cdot}\, d\mu<\infty,   \int_{X}\norm{\cdot}\, d\nu<\infty
\end{equation*}
and $\mu\prec_{\mathcal{F}}\nu$.
Suppose that any $\mathcal{F}$-transport between $\mu$ and $\nu$ is local.
Then there exists a Borel measurable set $B\subset X^*$ with $\mu(B)=1$ such that whenever $x_1^*,x_2^*\in B$ then 
\begin{equation*}
\mathrm{irc}_{\mathcal{A}}(\mu,\nu)(x_1^*)\cap\mathrm{irc}_{\mathcal{A}}(\mu,\nu)(x_2^*)=\emptyset
\end{equation*}
or
\begin{equation*}
\mathrm{irc}_{\mathcal{A}}(\mu,\nu)(x_1^*)=\mathrm{irc}_{\mathcal{A}}(\mu,\nu)(x_2^*).
\end{equation*}
These sets are relatively open, convex and finite-dimensional subsets of $X^*$.
Moreover, for any disintegration $\lambda\in\Lambda_{\mathcal{F}}(\mu,\nu)$
\begin{equation*}
    \mathrm{supp}\lambda(\omega,\cdot)\subset \mathrm{cl}\big(\mathrm{irc}_{\mathcal{A}}(\mu,\nu)(x^*)\big)\text{ for }\mu\text{-almost every }x^*\in X^*.
\end{equation*}
The sets $\mathrm{irc}_{\mathcal{A}}(\mu,\nu)(\cdot)$ are the smallest convex sets that satisfy the above condition.
Furthermore, if $\mathcal{B}=\mathrm{Cone}\mathcal{G}\cap(-\mathcal{G})$, then
\begin{equation*}
       \Phi_{\mathcal{B}}(\omega)\in R_{\mathcal{B}}\big(\mathrm{irc}_{\mathcal{A}}(\mu,\nu)(\omega)\big)\text{ for }\mu\text{-almost every }\omega\in\Omega.
        \end{equation*}
    
Let $U\subset\Omega\times\Omega$ be a Borel set such that for $\mu $-almost every $x^*\in X^*$
\begin{equation*}
U_{x^*}\subset\mathrm{irc}_{\mathcal{A}}(\mu,\nu)(x^*).
\end{equation*}
Then $U$ is $\Gamma_{\mathcal{F}}(\mu,\nu)$-polar set if and only if there exist Borel sets $N_1,N_2\subset\Omega$ with
\begin{equation*}
\mu(N_1)=0,\nu(N_2)=0
\end{equation*}
and 
\begin{equation*}
U\subset (N_1\times\Omega)\cup (\Omega\times N_2).
\end{equation*} 
\end{theorem}
\begin{proof}
     The theorem is a direct application of Theorem \ref{thm:partitionINTRO}, Theorem \ref{thm:polarinto}, Theorem \ref{thm:finestrucfinintro}. The verification of the assumptions of these theorems and the proof follows analogous lines to the proof of Theorem \ref{thm:martingale}.


\end{proof}

Measure $\mu$ may be disintegrated, thanks to \cite[Theorem 13.16]{Ciosmak20242}, in such a way that any conditional measure admits a single irreducible component.

If the considered transports are not assumed to be local, then we may still provide a partitioning, thanks to Section \ref{s:approx}, Theorem \ref{thm:apir}. We do not however obtain a characterisation of polar sets.

As in Section \ref{s:martingale} we may also provide a description of the intersections of the irreducible components thanks to Theorem \ref{thm:finestrucfinintro}.

\begin{theorem}\label{thm:intersectionssub}
   Suppose that $X^*$ is separable Banach space, dual to $X$. Let $\mu,\nu\in\mathcal{P}(\Omega)$ be two Radon measures such that 
\begin{equation*}
    \int_{X^*}\norm{\cdot}\, d\mu<\infty,   \int_{X}\norm{\cdot}\, d\nu<\infty
\end{equation*}
and $\mu\prec_{\mathcal{F}}\nu$.
Suppose that any $\mathcal{F}$-transport between $\mu$ and $\nu$ is local.
Then there exists a Borel measurable set $B\subset X^*$ with $\mu(B)=1$ such that whenever $x_1^*,x_2^*\in B$, and if 
\begin{equation*}
    x^*\in \mathrm{cl}\big(\mathrm{irc}_{\mathcal{A}}(\mu,\nu)(x_1^*)\big)\cap\mathrm{cl}\big(\mathrm{irc}_{\mathcal{A}}(\mu,\nu)(x_2^*)\big)
\end{equation*}
is  such that for $i=1,2$ there exists $\eta_i\in\mathcal{P}_{\norm{\cdot}}(\Omega)$ and $C_i>0$ such that $\eta_i\leq C_i\lambda(\omega_i,\cdot)$, where $\lambda\in\Lambda_{\mathcal{F}}(\mu,\nu)$ is a maximal disintegration, and
\begin{equation}\label{eqn:intrep}
x^*(x)=\int_{\Omega}x\, d\eta_i\text{ for all }x\in\mathcal{X},
\end{equation}
then we have the equality of the Gleason parts
\begin{equation*}
G\big(x^*,\mathrm{cl}\big(\mathrm{irc}_{\mathcal{A}}(\mu,\nu)(x_1^*)\big)\big)=G\big(x^*,\mathrm{cl}\big(\mathrm{irc}_{\mathcal{A}}(\mu,\nu)(x_2^*)\big)\big).
\end{equation*}
\end{theorem}
\begin{proof}
    The theorem is a consequence of Theorem \ref{thm:finestrucfinintro}.
\end{proof}

\begin{remark}
It is not necessarily true that $x^*\in \mathrm{irc}_{\mathcal{A}}(\mu,\nu)(x^*)$ for $\mu$-almost every $x^*\in X^*$. What is true, thanks to Remark \ref{rem:submart}, is that
\begin{equation*}
    x^*(x)\in \Phi(x)(\mathrm{irc}_{\mathcal{A}}(\mu,\nu)(x^*))\text{ for } x\in \mathrm{Cone}\mathcal{G}_0\cap (-\mathrm{Cone}\mathcal{G}_0)\text{ and }\mu\text{-almost every }x^*\in X^*.
\end{equation*}
\end{remark}

\begin{example}\label{exa:notcon}
    Let $\mu=\delta_{0}$ and let $\nu=\lambda|_{[1,2]}$, where $\lambda$ is the Lebesgue measure on $\mathbb{R}$. Let $\mathcal{G}$ consist of an increasing affine function. Then, as any $f\in\mathcal{F}$ is non-decreasing,
    \begin{equation*}
        \mu\prec_{\mathcal{F}}\nu.
    \end{equation*}
    Clearly $\mathrm{irc}_{\mathcal{A}}(\mu,\nu)(0)=(1,2)$. We see that $0\notin \mathrm{irc}_{\mathcal{A}}(\mu,\nu)(0)$.
\end{example}


If the considered transports are not necessarily local, then we may still provide a partitioning, thanks to Section \ref{s:approx}, Theorem \ref{thm:apir}. We do not however obtain a characterisation of polar sets.

\section{Approximating $\mathcal{F}$-transports}\label{s:relaxed}

In \cite[Definition 2.1., p. 4]{Guo2019}, the authors, given $\delta> 0$, define $\delta$-approximating martingale measure, which is an extension of the notion of a martingale measure. In this spirit we give the following definition.

\begin{definition}\label{def:delta}
    Suppose that $\mu,\nu\in\mathcal{P}_p(\Omega)$ are two Radon probability measures.
Let $\delta\geq 0$. 
Any Radon probability measure $\pi$ on $\Omega\times\Omega$ with marginals $\mu,\nu$ and such that for all $g\in\mathcal{G}$ and all bounded, non-negative measurable $h$ we have
\begin{equation}\label{eqn:delta}
\int_{\Omega}g(\omega_1)h(\omega_1)\,d\mu(\omega_1)-\int_{\Omega\times\Omega}g(\omega_2)h(\omega_1)\, d\pi(\omega_1,\omega_2)\leq \delta\norm{h}_{\mathcal{D}_1(\Omega)}\norm{g}_{\mathcal{D}_{p+1}(\Omega)}
\end{equation}
we shall call a $\delta$-approximating $\mathcal{F}$-transport between $\mu$ and $\nu$. 
We shall denote the set of all such measures by $\Gamma_{\mathcal{F}, \delta}(\mu,\nu)$. 
\end{definition}

\begin{remark}
    Clearly, $\Gamma_{\mathcal{F},\delta}(\mu,\nu)\supset\Gamma_{\mathcal{F}}(\mu,\nu)$. Therefore, if $\mu\prec_{\mathcal{F}}\nu$, then for all $\delta\geq 0$, $\Gamma_{\mathcal{F},\delta}(\mu,\nu)\neq\emptyset$. Note also that $\Gamma_{\mathcal{F},0}(\mu,\nu)=\Gamma_{\mathcal{F}}(\mu,\nu)$.
\end{remark}

\begin{remark}\label{rem:delta}
    As in \cite{Guo2019}, we observe that (\ref{eqn:delta}) is equivalent to assuming that for any $g\in\mathcal{G}$ we have
    \begin{equation*}
        \int_{\Omega}\Big(g(\omega_1)-\int_{\Omega}g\, d\lambda(\omega_1,\cdot)\Big)_+\, d\mu(\omega_1)\leq \delta\norm{g}_{\mathcal{D}_p(\omega)},
    \end{equation*}
    where $\lambda\in T(\pi)$ is a disintegration of $\pi$.
    If $\mathcal{G}$ is a linear subspace, then we note that (\ref{eqn:delta}) implies that
    \begin{equation*}
        \int_{\Omega}\Big\lvert g(\omega_1)-\int_{\Omega}g\, d\lambda(\omega_1,\cdot)\Big\rvert\, d\mu(\omega_1)\leq 2\delta\norm{g}_{\mathcal{D}_p(\omega)},
    \end{equation*}
    which recovers the setting of \cite{Guo2019} in the case of two measures in convex order.
\end{remark}


We shall now provide a precise formulation of the result concerning the finest partitioning into $\mathcal{A}$-convex sets in the setting of relaxed martingale transport, or, more generally, in the setting of $\delta$-approximating $\mathcal{F}$-transport.

\begin{definition}\label{def:deltap}
The set disintegration with respect to the first co-ordinate we shall denote by  $\Lambda_{\mathcal{F},\delta}(\mu,\nu)$. We shall also speak of maximal disintegration, which is an element $\lambda\in\Lambda_{\mathcal{F},\delta}(\mu,\nu)$ such that for all $\lambda'\in\Lambda_{\mathcal{F},\delta}(\mu,\nu)$ and for $\mu$-almost every $\omega\in\Omega$, 
\begin{equation*}
    \mathrm{supp}\lambda'(\omega,\cdot)\subset\mathrm{supp}\lambda(\omega,\cdot).
\end{equation*}
\end{definition}

We assume that $\Omega,\mathcal{G},\mathcal{F},\mathcal{A},p$ satisfy Assumption \ref{as:ass}. 

One may ask whether the irreducible paving can be also obtained in the setting of $\delta$-approximating $\mathcal{F}$-transports, and what the irreducible components will be. It can be shown that all the reasonings that worked for $\mathcal{F}$-transports will work as well for $\delta$-approximating $\mathcal{F}$-transports. The irreducible components were defined, in the case when $\mathcal{A}$ was finite-dimensional, as the relative interior of the closed convex hull of $\Phi(\mathrm{supp}\lambda(\omega,\cdot))$, where $\lambda\in\Lambda_{\mathcal{F}}(\mu,\nu)$ is a maximal disintegration. The following theorem shows however that in the case of $\delta$-approximating $\mathcal{F}$-transports the supports of any maximal disintegration are trivial.


\begin{theorem}\label{thm:relaxed}
Let $\delta> 0$.Suppose that $\mu,\nu$ are two Radon probability measure such that $\Lambda_{\mathcal{F},\epsilon}(\mu,\nu)$ is non-empty. Then any maximal disintegration $\lambda\in\Lambda_{\mathcal{F},\epsilon}(\mu,\nu)$ satisfies 
\begin{equation*}
    \mathrm{supp}\lambda(\omega,\cdot)=\mathrm{supp}\nu\text{ for }\mu\text{-almost every }\omega\in\Omega.
\end{equation*}
Moreover, a Borel set $U\subset\Omega\times\Omega$ is a polar set with respect to $\Gamma_{\mathcal{F},\delta}(\mu,\nu)$ if and only if there exist Borel sets $N_1,N_2\subset\Omega$ such that $\mu(N_1)=0$ and $\nu(N_2)=0$ and 
\begin{equation*}
    U\subset (N_1\times\Omega)\cup(N_2\times\Omega).
\end{equation*}
\end{theorem}
\begin{proof}
Let $\lambda\in\Lambda_{\mathcal{F}}(\mu,\nu)$ be any disintegration of an $\mathcal{F}$-transport. 
Consider for $t\in (0,1)$ the disintegration
\begin{equation*}
    \lambda_{t}=(1-t)\lambda+t\nu.
\end{equation*}
Then $\lambda\in \Lambda(\mu,\nu)$. Moreover for any $g\in\mathcal{G}$, 
\begin{equation*}
\int_{\Omega}\Big( g(\omega)-\int_{\Omega}g\, d\lambda_{t}(\omega,\cdot)\Big)_+\, d\mu(\omega)\leq t \int_{\Omega}\Big( g(\omega)-\int_{\Omega}g\, d\nu\Big)_+\, d\mu(\omega),
\end{equation*}
which can be bounded by 
\begin{equation*}
 t\int_{\Omega}\int_{\Omega}\abs{g(\omega_1)-g(\omega_2)}\, d\mu(\omega_1)\,d\nu(\omega_2)\leq t\norm{g}_{\mathcal{D}_p(\Omega)}\int_{\Omega}p\, d(\mu+\nu).
\end{equation*}
Thus if 
\begin{equation*}
    t\int_{\Omega}p\, d(\mu+\nu)\leq \epsilon,
\end{equation*}
then $\lambda_{t}$ is an $\delta$-approximating $\mathcal{F}$-transport, thanks to Remark \ref{rem:delta}. Its supports are equal to the support of $\nu$, thus maximal.

An analogous argument can be applied for the proof of the second assertion. Again, let $\lambda\in\Lambda_{\mathcal{F}}(\mu,\nu)$ be any disintegration. Note that any transport $\pi\in\Gamma(\mu,\nu)$ between $\mu$ and $\nu$, has finite cost 
\begin{equation}
    c(\pi)=\int_{\Omega\times\Omega}\sup\{(g(\omega_1)-g(\omega_2))_+\mid \norm{g}_{\mathcal{D}_p(\Omega)}\leq 1\} \, d\pi(\omega_1,\omega_2),
\end{equation}
as it is bounded by $\int_{\Omega}p\, d(\mu+\nu)$.
Let $\phi\in\Lambda(\mu,\nu)$ be a disintegration of $\pi$ with respect to the first co-ordinate.
Then for $t\in (0,1)$, considering $\lambda_t=(1-t)\lambda+t\phi$, we see that for any $g\in\mathcal{G}$ 
\begin{equation*}
\int_{\Omega}\Big( g(\omega)-\int_{\Omega}g\, d\lambda_{t}(\omega,\cdot)\Big)_+\, d\mu(\omega)
\end{equation*}
can be bounded by
\begin{equation*}
t \int_{\Omega}\Big( g(\omega)-\int_{\Omega}g\, d\phi\Big)_+\, d\mu(\omega)\leq t\int_{\Omega}\int_{\Omega}(g(\omega_1)-g(\omega_2))_+\,d\pi(\omega_1,\omega_2)\leq t\norm{g}_{\mathcal{D}_p(\Omega)}c(\pi).
\end{equation*}
It $t$ is sufficiently small, we see that $\lambda_t$ belongs to $\Lambda_{\mathcal{F},\delta}(\mu,\nu)$. Therefore, if $U$ is a polar set with respect to $\Gamma_{\mathcal{F},\delta}(\mu,\nu)$, then it is also polar with respect to $\Gamma(\mu,\nu)$. Then a classical result of Kellerer \cite[Proposition 3.3., p. 424 and Proposition 3.5., p. 425]{Kellerer1984}, also  see \cite[Proposition 2.1.]{Beiglbock2009},  completes the proof.
\end{proof}

This gives, in particular, a negative resolution of a question of Cox, \cite{Cox2022}, concerning the existence of such irreducible convex paving for the relaxed martingale transport on finite-dimensional linear spaces.

\section{Finite-dimensional approximations and general infinite-dimensional setting}\label{s:approx}

In Theorem \ref{thm:partitionINTRO} we assume that all $\mathcal{F}$-transports between $\mu$ and $\nu$ are local. Therefore, we do not immediately obtain a partitioning if there exists a non-local $\mathcal{F}$-transport between given two measures in $\mathcal{F}$-order.

A resolution of that obstacle can be achieved by the procedure of finite-dimensional approximations.

For $Z\subset\mathcal{G}$ we set $\mathcal{A}_Z$ to be the linear span of $Z$ and $\mathcal{F}_Z$ to be the lattice cone generated by $Z$. 
We define $\Phi_{\mathcal{A}_Z}\colon \Omega\to\mathcal{A}_Z$ via the formula
\begin{equation}
    \Phi_{\mathcal{A}_Z}(\omega)(a)=a(\omega)\text{ for }\omega\in\Omega\text{ and }a\in\mathcal{A}_Z.
\end{equation}
Since $\mathcal{F}_Z\subset\mathcal{F}$, if $\mu\prec_{\mathcal{F}}\nu$, then $\mu\prec_{\mathcal{F}_Z}\nu$ for any $Z\subset\mathcal{G}$. Therefore if $\mu\prec_{\mathcal{F}}\nu$, then for any $Z\subset\mathcal{G}$ and $\epsilon\geq 0$, the sets 
\begin{equation*}
\mathrm{irc}_{\mathcal{A}_Z}(\mu,\nu)(\cdot)
\end{equation*}
are well-defined subsets of $\mathcal{A}_Z^*$. Let us note that we do not need the assumption that $\mathcal{A}_Z$ separates points of $\Omega$, as this was used merely to infer the existence of transports in the proof of Theorem \ref{thm:strassenINTRO}, cf. \cite[Theorem 4.2]{Ciosmak20242}. By \cite[Theorem 13.6]{Ciosmak20242}, $\mathrm{irc}_{\mathcal{A}_Z}(\mu,\nu)(\cdot)$ consists of functionals of the form
\begin{equation*}
    \mathcal{A}_Z\ni a\mapsto\int_{\Omega}a\, d\eta\in\mathbb{R},
\end{equation*}
where $\lambda\in\Lambda_{\mathcal{F}_Z,\epsilon}(\mu,\nu)$ is a maximal disintegration,  $\eta\in \Xi(\lambda,\omega)$ is such that for some $c>0$, 
\begin{equation*}
    \eta\geq c\lambda(\omega,\cdot).
\end{equation*}
Each such functional admits a continuous extension to $\mathcal{A}$, thanks to the Hahn--Banach theorem. The set of all such extensions of functionals in $\mathrm{irc}_{\mathcal{A}_Z}(\mu,\nu)(\cdot)$ shall be denoted by $\mathrm{eirc}_{\mathcal{A}_Z}(\mu,\nu)(\cdot)$.

\begin{definition}\label{def:approximate}
      Suppose that $\mu\prec_{\mathcal{F}}\nu$. We define the approximate irreducible components $\mathrm{apirc}_{\mathcal{A}}(\mu,\nu)(\cdot)\subset\mathcal{A}^*$ by the formula
      \begin{equation*}
      \mathrm{apirc}_{\mathcal{A}}(\mu,\nu)(\omega)=\bigcap\Big\{\mathrm{eirc}_{\mathcal{A}_Z}(\mu,\nu)(\omega)\mid Z\subset\mathcal{G}\text{ is finite}\Big\}, \omega\in\Omega.
      \end{equation*}
\end{definition}
Below we shall use the restriction map $R_{\mathcal{B}}\colon \mathcal{A}^*\to \mathcal{B}^*$,
which assigns to a functional $a^*\in\mathcal{A}^*$ its restriction to a subspace $\mathcal{B}\subset\mathcal{A}$.
Moreover, $\xi\colon [0,\infty)\to\mathbb{R}$ is a positive function, of superlinear growth, such that $\xi(p)\in\mathcal{H}$ and $\xi(p)$ is integrable with respect to $\mu,\nu$. Existence of such $\xi$ follows by \cite[Theorem 2.10]{Ciosmak20242}.

\begin{lemma}\label{lem:closure}
Suppose that $\mu\prec_{\mathcal{F}}\nu$. Suppose that any $\mathcal{F}$-transport between $\mu$  and $\nu$ is local. Let $\lambda\in\Lambda_{\mathcal{F}}(\mu,\nu)$. Then for any measure $\eta\in\Xi(\lambda,\omega)$ such that for some $c>0$, 
\begin{equation*}
    \eta\geq c\lambda(\omega,\cdot),
\end{equation*}
the functional
\begin{equation*}
    \mathcal{A}\ni a\mapsto\int_{\Omega}a \,d\eta\in\mathbb{R}
\end{equation*}
belongs to the closure of the irreducible component $\mathrm{irc}_{\mathcal{A}}(\mu,\nu)$.
\end{lemma}
\begin{proof}
    Let $\lambda'\in\Lambda_{\mathcal{F}}(\mu,\nu)$ be a maximal disintegration of $\Gamma_{\mathcal{F}}(\mu,\nu)$. For $\epsilon\in (0,1)$, the disintegration $(1-\epsilon)\lambda+\epsilon\lambda'$ is also a maximal disintegration of $\Gamma_{\mathcal{F}}(\mu,\nu)$. Letting $\epsilon$ tend to zero completes the proof.
    \end{proof}

Below we shall state a theorem that yields a partitioning in the infinite-dimensional setting. We would like for the elements of the partitioning to be the approximate irreducible components, see Definition \ref{def:approximate}. We do not know how to prove that these are pairwise disjoint on a set of full measure $\mu$ -- for the proof of this we would need e.g. semi-continuity of the sets $\mathrm{eirc_{\mathcal{A}_Z}}(\mu,\nu)(\cdot)$ with respect to $Z$.
Instead, we give a proof of an analogous theorem for a family of larger sets
 \begin{equation}\label{eqn:apir}
      \mathrm{apirc}_{\mathcal{G}_0}(\mu,\nu)(\omega)=\bigcap\Big\{\mathrm{eirc}_{\mathcal{A}_Z}(\mu,\nu)(\omega)\mid Z\subset\mathcal{G}_0\text{ is finite}\Big\}, \omega\in\Omega,
      \end{equation}
      where $\mathcal{G}_0$ is at most countable subset of $\mathcal{G}$.

\begin{theorem}\label{thm:apir}
    Suppose that $\mu\prec_{\mathcal{F}}\nu$. Suppose that $\mathcal{G}_0$ is at most countable subset of $\mathcal{G}$. Then there exists a Borel measurable set $B\subset\Omega$ with $\mu(B)=1$ such that whenever $\omega_1,\omega_2\in B$  then
    \begin{equation*}
        \mathrm{apirc}_{\mathcal{G}_0}(\mu,\nu)(\omega_1)\cap \mathrm{apirc}_{\mathcal{G}_0}(\mu,\nu)(\omega_2)=\emptyset
    \end{equation*}
    or
    \begin{equation*}
        \mathrm{apirc}_{\mathcal{G}_0}(\mu,\nu)(\omega_1)= \mathrm{apirc}_{\mathcal{G}_0}(\mu,\nu)(\omega_2).
    \end{equation*}
    Moreover, the sets $\mathrm{apirc}_{\mathcal{G}_0}(\mu,\nu)(\cdot)$ are convex and for any  $\lambda\in\Lambda_{\mathcal{F}}(\mu,\nu)$
    \begin{equation}\label{eqn:constr}
        \mathrm{supp}\lambda(\omega,\cdot)\subset \Phi^{-1}\Big(\mathrm{cl}\big(\mathrm{apirc}_{\mathcal{G}_0}(\mu,\nu)(\omega)\big)\Big)\text{  for }\mu\text{-almost every }\omega\in\Omega.
    \end{equation}
Furthermore, for $\mathcal{B}=\mathrm{Cone}\mathcal{G}\cap(-\mathrm{Cone}\mathcal{G})$,
    \begin{equation*}
       \Phi_{\mathcal{B}}(\omega)\in R_{\mathcal{B}}\big(\mathrm{cl}\big(\mathrm{apirc}_{\mathcal{G}_0}(\mu,\nu)(\omega)\big)\big)\text{ for }\mu\text{-almost every }\omega\in\Omega.
        \end{equation*}
    In particular, if $\mathcal{G}$ is symmetric, then 
    \begin{equation*}
        \Phi(\omega)\in \mathrm{cl}\big(\mathrm{apirc}_{\mathcal{G}_0}(\mu,\nu)(\omega)\big)\text{ for }\mu\text{-almost every }\omega\in\Omega.
    \end{equation*}
\end{theorem}

\begin{proof}
For any finite $Z\subset\mathcal{G}_0$ the space $\mathcal{A}_Z$ is finite-dimensional.  Then \cite[Theorem 13.6]{Ciosmak20242} shows that the set 
\begin{equation*}
    \mathrm{irc}_{\mathcal{A}_Z}(\mu,\nu)(\omega)
\end{equation*}
consists of all functionals 
\begin{equation*}
    \mathcal{A}_Z\ni a\mapsto\int_{\Omega}a\, d\eta\in\mathbb{R},
\end{equation*}
for any measure $\eta\in\Xi(\lambda,\omega)$ such that for some $c>0$, 
\begin{equation*}
    \eta\geq c\lambda(\omega,\cdot),
    \end{equation*}
 and $\lambda\in\Lambda_{\mathcal{F}_Z}(\mu,\nu)$ is a maximal disintegration of $\Gamma_{\mathcal{F}_Z}(\mu,\nu)$, the set of $\mathcal{F}_Z$-transports between $\mu$ and $\nu$. Therefore $\mathrm{eirc}_{\mathcal{A}_Z}(\mu,\nu)(\omega)$ contains the functionals
\begin{equation*}
    \mathcal{A}\ni a\mapsto\int_{\Omega}a\, d\eta\in\mathbb{R},
\end{equation*}
for $\eta$ as above.

Since $\Lambda_{\mathcal{F}}(\mu,\nu)\subset\Lambda_{\mathcal{F}_Z}(\mu,\nu)$ for any $\lambda\in \Lambda_{\mathcal{F}}(\mu,\nu)$, $Z\subset\mathcal{G}_0$, the elements $a^*\in\mathcal{A}^*$ such that 
\begin{equation*}
a^*(a)=\int_{\Omega}a \,d\eta\text{ for all }    a\in \mathcal{A}, 
\end{equation*}
for some Radon probability measure $\eta\in\mathcal{P}_{\xi(p)}(\Omega)$ such that for some $C,c>0$, 
\begin{equation*}
 C\lambda(\omega,\cdot)\geq   \eta\geq c\lambda(\omega,\cdot),
\end{equation*}
belong to $\mathrm{apirc}_{\mathcal{A},0}(\mu,\nu)(\omega)$.
These functionals are dense in $\Phi(\mathrm{supp}\lambda(\omega,\cdot))$, see \cite[Lemma 13.5]{Ciosmak20242}. This shows that (\ref{eqn:constr}) holds true.


We see that in the definition of $\mathrm{apirc}_{\mathcal{G}_0}(\mu,\nu)(\cdot)$ we take the intersections over a countable set $Y$ of finite sets $Z\subset\mathcal{G}_0$. For each such $Z\subset\mathcal{G}_0$ we may pick a Borel measurable set $B_{Z}$ with $\mu(B_{Z})=1$ and such that for $\omega_1,\omega_2\in B_{Z}$ the irreducible components $\mathrm{irc}_{\mathcal{A}_Z}(\mu,\nu)(\omega_i)$, $i=1,2$, are either equal or disjoint. Let $B$ be the intersection over the countable family of sets $(B_{Z})_{Z\in Y}$. Then whenever $\omega_1,\omega_2\in B$, then the approximate irreducible components $\mathrm{apirc}_{\mathcal{G}_0}(\mu,\nu)(\omega_i)$, $i=1,2$ are either equal or disjoint. Moreover, $\mu(B)=1$. 

Indeed, we need to show that for $\omega_1,\omega_2\in B_{Z}$, the sets $\mathrm{eirc}_{\mathcal{A}_Z}(\mu,\nu)(\omega_i)$, $i=1,2$, are either equal or disjoint.
If these sets intersect at some $a^*\in\mathcal{A}^*$, then the restriction of $a^*$ to $\mathcal{A}_Z$, belongs to the both irreducible components $\mathrm{irc}_{\mathcal{A}_Z}(\mu,\nu)(\omega_i)$, $i=1,2$. By the property of $B_{Z}$, this implies that the irreducible components are equal, and therefore so are the sets $\mathrm{eirc}_{\mathcal{A}_Z}(\mu,\nu)(\omega_i)$, $i=1,2$.

    If $g\in\mathcal{G}\cap(-\mathcal{G})$ and $\lambda\in\Lambda_{\mathcal{F}}(\mu,\nu)$, 
    then for $\mu$-almost every $\omega\in\Omega$ 
    \begin{equation*}
         g(\omega)=\int_{\Omega}g\, d\lambda(\omega,\cdot).
    \end{equation*}
      The functional 
    \begin{equation*}
        \mathcal{A}\ni a\mapsto \int_{\Omega}a\, d\lambda(\omega,\cdot)\in\mathbb{R}
    \end{equation*}
    belongs to the closure of $\mathrm{apirc}_{\mathcal{G}_0}(\mu,\nu)$, by Lemma \ref{lem:closure}.
    This shows that 
    \begin{equation*}
        \Phi_{\mathcal{B}}(\omega)\in R_{\mathcal{B}}\big(\mathrm{cl}\big(\mathrm{apirc}_{\mathcal{A},0}(\mu,\nu)(\omega)\big)\big)\text{ for }\mu\text{-almost every }\omega\in\Omega.
    \end{equation*}
    The last assertion is a consequence of the fact that if $\mathcal{G}$ is symmetric then $\mathcal{B}=\mathcal{A}$.
    
    
    
\end{proof}


\begin{remark}
    Theorem \ref{thm:apir} allows to extend the results of Section \ref{s:martingale} and Section \ref{s:supermartingale} to the non-local setting, in an obvious way. We do not however obtain a characterisation of polar sets.
\end{remark}

\begin{remark}\label{rem:dense}
Suppose that the countable subset $\mathcal{G}_0\subset\mathcal{G}$ is dense in $\mathcal{G}$ in the norm of $\mathcal{D}_{p+1}(\Omega)$. We conjecture that under this assumption the approximate irreducible components defined in Definition \ref{def:approximate} and in (\ref{eqn:apir}) coincide.
\end{remark}

\section{Harmonic transports and Skorokhod embedding problem}\label{s:harmonic}

Let us now consider transports related to the Skorokhod embedding problem. 

\begin{definition}
Let $\Omega\subset\mathbb{R}^n$ be an open set. We shall say that a continuous function $a$ on $\Omega$ is harmonic whenever the value of $a$ at any point is equal to the average of $a$ on any  ball contained in $\Omega$ and centred at the point, i.e., for all $\omega\in\Omega$ and $r>0$ such that $B(\omega,r)\subset\Omega$
\begin{equation*}
    a(\omega)= \frac1{\lambda(B(\omega,r))}\int_{B(\omega,r)}a\, d\lambda,
\end{equation*}
where $B(\omega,r)\subset\Omega$ is a ball of radius $r$ centred at $\omega\in\Omega$, and $\lambda$ is the Lebesgue measure.
Let $\mathcal{A}$ denote the linear space of continuous harmonic functions on $\Omega$.
An upper semi-continuous continuous function is said to be subharmonic whenever the value at any point is bounded from the above by the average of $h$ on any  ball contained in $\Omega$ and centred at the point, i.e., for all $\omega\in\Omega$ and $r>0$ such that $B(\omega,r)\subset\Omega$
\begin{equation}\label{eqn:upper}
    h(\omega)\leq \frac1{\lambda(B(\omega,r))}\int_{B(\omega,r)}h\, d\lambda.
\end{equation}
Let $\mathcal{S}$ denote the cone of upper semi-continuous subharmonic functions on an open set  $\Omega\subset\mathbb{R}^n$. 
\end{definition}

Given two Radon probability measures $\mu,\nu$ we say that they are in subharmonic order, whenever $\mu\prec_{\mathcal{S}}\nu$.

Let us note the relation of measures in subharmonic order to the Skorokhod embedding problem, see e.g. \cite{Palmer2019}.
It tells that $\mu\prec_{\mathcal{S}}\nu$ if and only if there exists a stopping time $\tau$ and a Brownian motion $(B_t)_{t\in [0,\infty)}$ such that $B_0\sim \mu$ and $B_{\tau}\sim \nu$. 
The existence of such a stopping time has been investigated by a great many studies, see e.g. \cite{Rost1971, Monroe1972, Baxter1974, Walsh1976, Falkner1980}. We refer the reader to the survey \cite{Obloj2004} of Ob\l\'oj for a thorough discussion of the Skorokhod embedding problem.

Let us mention a generalisation of the subharmonic transport problem, discussed in \cite{Kim2021}, where the authors consider general Feller processes in place of Brownian motion.
Our theory is flexible enough to provide a partitioning also in that setting.

Let us note that clearly that the cone of continuous functions in the complete lattice cone $\mathcal{H}$ generated by $\mathcal{A}$ is contained in the cone $\mathcal{R}$ of lower semi-continuous functions that satisfy (\ref{eqn:upper}). Moreover $\mathcal{S}\cap (-\mathcal{S})=\mathcal{A}$, by the very definition. 
However, $\mathcal{H}$ and and the cone $\mathcal{R}$ do not coincide. Indeed, if they did, then any Radon probability measure $\eta$ on $\Omega$ such that for some $\omega\in\Omega$ 
\begin{equation*}
    a(\omega)=\int_{\Omega}a\, d\eta\text{ for all }a\in\mathcal{A},
\end{equation*}
would satisfy 
\begin{equation*}
    r(\omega)\leq\int_{\Omega}r\, d\eta\text{ for all }r\in\mathcal{R}.
\end{equation*}
This is not true in general, as shown by \cite[Example 5.2.7., p. 83]{Alan2009}.

There are therefore two possible approaches that we may follow. One is to consider irreducible pavings for pairs of measures in subharmonic order. The other is to consider irreducible pavings for pairs of measures in order with respect to the complete lattice cone generated by the space of harmonic functions.
Despite it is less relevant to the study of the Skorokhod embedding problem, we shall follow the latter path, leaving the former one for the future research.

Let us observe that Theorem \ref{thm:holoin} shows that any open set $\Omega\subset\mathbb{R}^n$ admits a continuous function $p\colon\Omega\to [0,\infty)$ that is proper and belongs to $\mathcal{H}$. Indeed, $\mathcal{A}$ is a closed subspace with respect to the compact-open topology on the space $\mathcal{C}(\Omega)$ of continuous functions on $\Omega$. Therefore existence of such $p$ is implied by the fact that $\Omega$ is complete with respect to $\mathcal{A}$. Consider the set of fundamental solution of the Laplace operator with poles at the boundary points of $\Omega$. They are continuous harmonic functions on $\Omega$, which verify that any Cauchy sequence with respect to $\mathcal{A}$ has its limit in $\Omega$. This is to say, $\Omega$ is complete with respect to $\mathcal{A}$.
Such a function $p$ clearly belongs also to $\mathcal{S}$. We may arrange that $p$ majorises any chosen finite family of functions in  $\mathcal{A}$ -- we assume that $p$ majorises affine functions.

\begin{theorem}\label{thm:harmonic}
Let $\Omega$ be an open subset of $\mathbb{R}^n$ and let $\mathcal{A}$ denote the space of continuous harmonic functions on  $\Omega$. Let $\mu,\nu$ be two Radon, probability measures on $\Omega$ in  order with respect to the lattice cone $\mathcal{F}$ generated by  $\mathcal{A}$, i.e.,
\begin{equation*}
\int_{\Omega}f\, d\mu\leq\int_{\Omega}f\, d\nu,
\end{equation*}
for all  functions $f\in\mathcal{F}$ integrable with respect to $\mu$ and $\nu$. 
Let $p\in\mathcal{H}$ be a non-negative, proper, continuous function on $\Omega$ in the complete lattice cone $\mathcal{H}$ generated by $\mathcal{A}$. Assume that 
\begin{equation*}
\int_{\Omega}p\, d\mu<\infty\text{ and }\int_{\Omega}p\, d\nu<\infty.
\end{equation*}
Let $\mathcal{B}$ denote the linear space of all functions in $\mathcal{A}$ that are of $(p+1)$-growth, let $\mathcal{E}$ denote the lattice cone generated  by $\mathcal{B}$.
Then the set of Radon measures $\pi$ on $\Omega\times\Omega$ with marginals $\mu$ and $\nu$ such that for all  $e\in\mathcal{E}$ of $(p+1)$-growth, and bounded, non-negative, measurable function $g$,
\begin{equation*}
\int_{\Omega}g(\omega_1)e(\omega_1)\, d\mu(\omega_1)\leq \int_{\Omega\times\Omega}g(\omega_1)e(\omega_2)\, d\pi(\omega_1,\omega_2)
\end{equation*}
is non-empty.
Moreover, to any $\omega\in\Omega$ we may assign a convex set $\mathrm{apirc}_{\mathcal{G}_0}(\mu,\nu)(\omega)\subset\mathcal{B}^*$, such that there exists a Borel set $B\subset\Omega$, $\mu(B)=1$ such that whenever $\omega_1,\omega_2\in B$ then
\begin{equation*}
    \mathrm{apirc}_{\mathcal{G}_0}(\mu,\nu)(\omega_1)\cap \mathrm{apirc}_{\mathcal{G}_0}(\mu,\nu)(\omega_2)\neq\emptyset
    \end{equation*}
implies that
    \begin{equation*}
    \mathrm{apirc}_{\mathcal{G}_0}(\mu,\nu)(\omega_1)=\mathrm{apirc}_{\mathcal{G}_0}(\mu,\nu)(\omega_2).
\end{equation*}
Furthermore, for any disintegration $\lambda\in\Lambda_{\mathcal{E}}(\mu,\nu)$
\begin{equation*}
    \mathrm{supp}\lambda(\omega,\cdot)\subset\Phi^{-1}\big((\mathrm{cl}\big(\mathrm{apirc}_{\mathcal{G}_0}(\mu,\nu)(\omega)\big)\big)\text{ for }\mu\text{-almost every }\omega\in\Omega,
\end{equation*}
and
\begin{equation*}
  \omega\subset\Phi^{-1}\big((\mathrm{cl}\big(\mathrm{apirc}_{\mathcal{G}_0}(\mu,\nu)(\omega)\big)\big)\text{ for }\mu\text{-almost every }\omega\in\Omega.
\end{equation*}
\end{theorem}
\begin{proof}
The first part of the theorem follows from \cite[Theorem 4.2]{Ciosmak20242}.
To verify its assumptions, it suffices to notice the Euclidean topology of $\Omega$ and the topology $\tau(\mathcal{B})$ coincide, since we assumed that $p$ majorises affine functions. 

The existence of the sets $\mathrm{apirc}_{\mathcal{B},0}(\mu,\nu)(\cdot)$ with the required properties follows from Theorem \ref{thm:apir}.
\end{proof}

\begin{remark}
    By \cite[Theorem 2.10]{Ciosmak20242} there exists $\xi$, increasing, positive, convex, of superlinear growth, such that $\mu$ and $\nu$ are $\xi(p)$-integrable and such that $\xi(p)\in\mathcal{H}$. 
Observe also that $\mathcal{B}$ is a separable subset of $\mathcal{D}_{\xi(p),0}(\Omega)$. Indeed, the space $\mathcal{C}_{\xi(p),0}(\Omega)$ itself is separable by the Stone--Weierstrass theorem and local compactness of $\Omega$, and so is any its subset. Therefore, we may take a countable, dense subset $\mathcal{G}_0\subset\mathcal{B}$ when applying above Theorem \ref{thm:apir}.
\end{remark}

\begin{example}\label{exa:harmonic}
Note that any two measures in subharmonic order are also in convex order. Therefore, to any such pair of measures, we may apply the results of Section \ref{s:martingale} and Theorem \ref{thm:harmonic}.
We shall provide an example of two Radon probability measures $\mu,\nu$ on $\mathbb{R}^2$ in subharmonic order such that the convex partitioning and the subharmonic partitioning differ. 

Let 
\begin{equation*}
    \mu=\frac12(\delta_{(0,0)}+\delta_{(0, 2)}),\, \nu=\frac12(\sigma_1+\sigma_2),
\end{equation*}
where $\sigma_1$ is the uniform measure on the circle $S$ of radius one, centred at $(0,0)$, and $\sigma_2$ is the distribution of $B_{\tau}$, where $(B_t)_{t\in [0,\infty)}$ is a standard Brownian motion started at $(0,2)$, $\tau$ is the first exit time of $(B_t)_{t\in [0,\infty)}$ from a smoothly bounded domain $U\subset\mathbb{R}^2$ that is disjoint from the disc of radius two centred at $(0,-1.5)$, but such that it contains the unit disc in the relative interior of its convex hull. The situation is illustrated below.

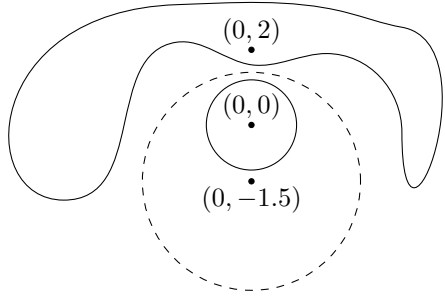
\begin{tikzpicture}
\coordinate (A) at (0,0);
\coordinate (C) at (0,-0.75);
\coordinate (B) at (0,1);
\draw  (0,0) circle (0.6);
\draw [dashed] (0,-0.75) circle (1.45);
\draw[fill] (A) circle (1.0pt);
\draw[fill] (C) circle (1.0pt);
\draw[fill] (B) circle (1.0pt);
\node[yshift=7pt] at (A) {$(0,0)$};
\node[yshift=-7pt] at (C) {$(0,-1.5)$};
\node[yshift=7pt] at (B) {$(0,2)$};
\coordinate (a1) at (2,-0.1);
\coordinate (a2) at (2,1.3);
\coordinate (a3) at (1.25,1.5);
\coordinate (a4) at (-1,1.6);
\coordinate (a5) at (-2.5,-1);
\coordinate (a6) at (-1,1.1);
\coordinate (a7) at (0,0.8);
\coordinate (a8) at (1,0.95);
\draw[use Hobby shortcut] ([out angle=-90]a1)..(a2)..(a3)..(a4)..(a5)..(a6)..(a7)..(a8)..([in angle=90]a1); 
\end{tikzpicture}

We claim that there is only one irreducible component pertaining to martingale transports between $\mu$ and $\nu$. It is the relative interior of the convex hull of $U$. Indeed, if
\begin{equation*}
    \lambda((0,0),\cdot)=\sigma_1,
    \lambda((0,2),\cdot)=\sigma_2
\end{equation*}
then $\lambda\in \Lambda_{\mathcal{S}}(\mu,\nu)$.
Let $\mathcal{C}$ denote the space of affine functions and let $\mathcal{D}$ be the corresponding lattice cone generated by $\mathcal{C}$. Then $\lambda\in \Lambda_{\mathcal{D}}(\mu,\nu)$ and therefore
\begin{equation*}
    \mathrm{irc}_{\mathcal{C}}(\mu,\nu)(0,2)=\mathrm{rintclConv}U,
\end{equation*}
and 
\begin{equation*}
    \mathrm{irc}_{\mathcal{C}}(\mu,\nu)(0,0)\supset \mathrm{rintclConv}S.
\end{equation*}
We see that the components intersect, which yields that they are equal, by Theorem \ref{thm:partitionINTRO}. The claim is proven.

Note however that there are two distinct approximate irreducible components pertaining to harmonic transports between $\mu$ and $\nu$. Indeed, let $D$ be the unit open disc of radius one, centred at the origin.
Let us consider 
\begin{equation*}
    h(\omega)= -\log \norm{\omega-\omega_0}\text{, for }\omega\in\mathbb{R}^2\setminus \{\omega_0\}\text{, where }\omega_0=(0,-1.5).
\end{equation*}
Then $h$ is harmonic on its domain. By the Runge theorem we may approximate it by harmonic polynomials on a compact set that contains $U\cup D$ in its interior. We may pick a harmonic polynomial $a$ in such a way that its level set that separates $D$ and $U$, as this was true for $h$. 
Therefore, for some $t\in\mathbb{R}$
\begin{equation*}
    D\subset a^{-1}(t,\infty)), U\subset a^{-1}((-\infty, t)).
\end{equation*}
Then $g=a\vee t$ is a continuous subharmonic function for which 
\begin{equation*}
    \int_{\mathbb{R}^n}g\, d\mu=\int_{\mathbb{R}^n}g\, d\nu.
\end{equation*}
If $\lambda\in\Lambda_{\mathcal{E}}(\mu,\nu)$, then
\begin{equation*}
   g((0,0))=\int_{\mathbb{R}^2}g\, d\lambda((0,0),\cdot) \text{ and }g((0,2))=\int_{\mathbb{R}^2}g\, d\lambda((0,2),\cdot).
\end{equation*}

Therefore, we see that $\lambda((0,0),\cdot)$ is supported on the set $g^{-1}((t,\infty))$ and $\lambda((0,2),\cdot)$ is supported on the set $g^{-1}((-\infty,t])$. This implies that
\begin{equation*}
    \mathrm{apirc}_{\mathcal{G}_0}(\mu,\nu)(0,0)\subset \mathrm{cl}D\text{ and } \mathrm{apirc}_{\mathcal{G}_0}(\mu,\nu)(0,2)\subset\mathrm{cl} U,
\end{equation*}
where $\mathcal{G}_0$ is a set containing $g$.
In particular, these sets are disjoint, so that our assertion is proven.
\end{example}



\section{Martingale problems}\label{s:martingaleproblem}

One of  the applications of the developments of \cite{Ciosmak20242} pertains to the martingale problems. We refer the reader to the book of Stroock and Varadhan \cite{Stroock2006} for an introduction to the theory of martingale problems. Let $\mathcal{B}(\Omega)$ denote the set of all Borel measurable functions on a set $\Omega\subset\mathbb{R}^n$, and let $\mathcal{D}\subset\mathcal{B}(\Omega)$ be a subspace. Let $L\colon \mathcal{D}\to\mathcal{B}(\Omega)$ be a linear operator and let $\mu\in\mathcal{P}(\Omega)$ be a Borel probability measure on $\Omega$. We shall say that a stochastic process $(X_t)_{t\in [0,1]}$ on a probability space equipped with a probability measure $\mathbb{P}$, with values in $\Omega$, solves the martingale problem for $(L,\mu)$ if $X_0\sim \mu$ and for any $f\in\mathcal{D}$
\begin{equation*}
    f(X_t)-\int_0^t (Lf)(X_s)\, ds
\end{equation*}
is a well-defined martingale with respect to the natural filtration $(\mathcal{F}_t)_{t\in [0,1]}$, that is 
\begin{equation*}
    \mathcal{F}_t=\sigma(\{X_s\mid s\in [0,t]\})\text{ for }t\in [0,1].
\end{equation*}

Suppose that $(X_t)_{t\in [0,1]}$ is a solution to the martingale problem for $(L,\mu)$. Let $\nu$ be the distribution of $X_1$. Let $\mathcal{A}$ be the kernel of $L$ and let $\mathcal{F}$ be the lattice cone generated  by $\mathcal{A}$. By the assumption for any $a\in \mathcal{A}$, $(a(X_t))_{t\in [0,1]}$ is a martingale. In particular, $(a(X_0),a(X_1))$ is a one-step martingale. By the results of Section \ref{s:approx}, Theorem \ref{thm:apir}, we see that for any countable subset $\mathcal{Z}\subset \mathcal{A}$ we may obtain a partitioning of our space into $\mathcal{A}$-convex sets $(\Phi^{-1}(\mathrm{apirc}_{\mathcal{Z}}(\mu,\nu)(\omega)))_{\omega\in\Omega}$, that yields constraints for the behaviour of the process $(X_0,X_1)$. That is, for any $\lambda\in\Lambda_{\mathcal{F}}(\mu,\nu)$,
\begin{equation*}
    \mathrm{supp}\lambda(\omega,\cdot)\subset\Phi^{-1}(\mathrm{cl}(\mathrm{apirc}_{\mathcal{Z}}(\mu,\nu)(\omega)))\text{ for }\mu\text{-almost every }\omega\in\Omega,
\end{equation*}
and 
\begin{equation*}
 \omega\in \Phi^{-1}(\mathrm{cl}(\mathrm{apirc}_{\mathcal{Z}}(\mu,\nu)(\omega)))\text{ for }\mu \text{-almost every }\omega\in\Omega.
\end{equation*}
In the probabilistic language the above translate to
\begin{equation*}
   X_1\in \Phi^{-1}(\mathrm{cl}(\mathrm{apirc}_{\mathcal{Z}}(\mu,\nu)(X_0)))\text{ almost surely,}
\end{equation*}
and
\begin{equation*}
 X_0\in \Phi^{-1}(\mathrm{cl}(\mathrm{apirc}_{\mathcal{Z}}(\mu,\nu)(X_0)))\text{ almost surely.}
\end{equation*}
The components give also constraints for the entire process $(X_t)_{t\in [0,1]}$. Indeed, for any $t\in [0,1]$ and $a\in\mathcal{A}$
\begin{equation*}
    a(X_t)=\mathbb{E}(a(X_1)\mid \mathcal{F}_t).
\end{equation*}
Thus for any  $t\in [0,1]$
\begin{equation}\label{eqn:t}
   X_t\in \Phi^{-1}(\mathrm{cl}(\mathrm{apirc}_{\mathcal{Z}}(\mu,\nu)(X_0)))\text{ almost surely.}
\end{equation}
Section \ref{s:harmonic}, in particular Example \ref{exa:harmonic}, shows that the partitioning and the obtained constraints can be non-trivial. 

Let us remark that one may, for any $t\in [0,1]$, consider the partitioning for the pair $(X_0,X_t)$. Then (\ref{eqn:t}) shows that the corresponding partitionings are increasing, with respect to the order induced by inclusion.

As we have mentioned in the introduction, there are no restrictions concerning the choice of subspace $\mathcal{A}$. If we take $\mathcal{A}$ to be the space of real parts of holomorphic functions, then we obtain a holomorphically convex paving for any holomorphic martingale. 

\section{Optimal transports}\label{s:optimal}


Let us demonstrate how the developed localisation scheme works for the case of optimal transports.

We shall consider an optimal transport problem with metric cost. Let $\Omega$ be a Riemannian manifold with metric $d$. Let $\mu,\nu$ be two Borel probability measures on $\Omega$, absolutely continuous with respect to the Lebesgue measure, with finite first moments, i.e.,
\begin{equation*}
    \int_{\Omega}d(\cdot,\omega_0)\, d\mu<\infty , \int_{\Omega}d(\cdot,\omega_0)\, d\nu<\infty 
\end{equation*}
for some, or equivalently any, $\omega_0\in\Omega$.

We shall show that our theory includes the localisation technique developed in \cite{Klartag2017} and further  in \cite{Cavalletti2017, Mondino2017, Ohta2018, Ciosmak2021, Ciosmak20212}.

By the Kantorovich  duality, see \cite{Kantorovich1958}, the optimal transport cost
\begin{equation*}
    \inf\Big\{\int_{\Omega\times\Omega}d(\omega_1,\omega_2)\, d\pi(\omega_1,\omega_2)\mid \pi\in\Gamma(\mu,\nu)\Big\}
\end{equation*}
coincides with  
\begin{equation*}
    \sup\Big\{\int_{\Omega}u\, d(\nu-\mu)\mid u\colon\Omega\to\mathbb{R}\text{ is }1\text{-Lipschitz}\Big\}.
\end{equation*}
It is well-known that the supremum is attained, see e.g. \cite[Theorem 2, p. 230]{Ciosmak20211}. A $1$-Lipschitz function $v$ that attains the supremum we shall call an optimal potential. 

Let $v\colon\Omega\to\mathbb{R}$ be an optimal potential for the transport between $\mu$ and $\nu$. Let us consider the tangent cone $\mathcal{F}$ at $1$-Lipschitz function $v$, to the set of all $1$-Lipschitz functions on $\Omega$, i.e.
\begin{equation*}
\mathcal{F}=\{\lambda(v-u)\mid \lambda\geq 0,u\text{ is }1\text{-Lipschitz}\}.
\end{equation*}
It is straightforward to verify that if $v$ is an optimal potential for $\mu$ and $\nu$, then $\mu\prec_{\mathcal{F}}\nu$, see e.g. \cite{Ciosmak20232}.

Thanks to the assumption on the finiteness of the first moments of $\mu$ and  $\nu$, we see that 
\begin{equation*}
    \mathcal{F}\subset L^1(\Omega,\mu+\nu).
\end{equation*}
It is immediate that $\mu\prec_{\mathrm{cl}\mathcal{F}}\nu$, where the closure is taken with respect to the metric of $L^1(\Omega,\mu+\nu)$.

We shall look for a linear subspace $\mathcal{A}\subset\mathrm{cl}\mathcal{F}$.

Let  us recall that for a $1$-Lipschitz function $u\colon\Omega\to\mathbb{R}$, a transport ray is a maximal set $\mathcal{T}\subset\Omega$ such that the restriction of $u$ to $\mathcal{T}$ is an isometry, cf. \cite[Definition 2.2.]{Ciosmak2021} and \cite[Definition 2.1.2.]{Klartag2017}, i.e., 
\begin{equation*}
    \abs{u(\omega_1)-u(\omega_2)}=d(\omega_1,\omega_2)\text{ for all  }\omega_1,\omega_2\in\mathcal{T}.
\end{equation*}

\begin{proposition}\label{pro:subspace}
    The maximal linear subspace $\mathcal{B}\subset\mathrm{cl}\mathcal{F}$ of Borel measurable functions in $L^1(\Omega,\mu+\nu)$ contains the space $\mathcal{A}$ all integrable, Borel measurable functions that are constant on transport rays of $v$.
    Moreover 
    \begin{equation*}
        \mathrm{cl}(\mathcal{F}\cap(-\mathcal{F}))\subset \mathcal{A}.
    \end{equation*}
\end{proposition}
\begin{proof}

Let us observe that $\mathcal{B}=\mathrm{cl}\mathcal{F}\cap (-\mathrm{cl}\mathcal{F})$.

    Clearly, the space $\mathcal{A}$ of all integrable, Borel measurable functions that are constant on transport rays of $v$ is a linear space. We shall show that it is contained in $\mathrm{cl}\mathcal{F}$. 
    We shall follow the strategy similar to the one in \cite[Lemma 4.3., Lemma 4.4., p. 56, 57]{Klartag2017}.

Note that \cite[Lemma 4.2.]{Klartag2017}, cf. \cite[Corollary 2.15., p. 224]{Ciosmak2021} for the Euclidean case, shows that the set of points that belong to at least two distinct transport rays have the Lebesgue measure zero. Let $A\subset \Omega$ be a Borel set that is a union of some transport rays of $v$. Let $\epsilon>0$, and let   $K^{\epsilon}\subset A$ be a compact set such that with $\mu(A\setminus K^{\epsilon})+\nu(A\setminus K^{\epsilon})\leq \epsilon$. Pick $\delta>0$ and set $v_{\delta,0}=v-\delta$ on $K^{\epsilon}$ and $v_{\delta,0}=v$ on 
\begin{equation*}
K_{\delta}^{\epsilon}=\{\omega\in \Omega\mid \delta\leq d(\omega,\omega')-\abs{v(\omega)-v(\omega')}\text{ for all }\omega'\in K^{\epsilon}\}.
\end{equation*} 
Then $v_{\delta,0}$ is $1$-Lipschitz on $K^{\epsilon}\cup K_{\delta}^{\epsilon}$ and within $\delta$-distance to $v$. We may extend it to $\Omega$ to $v_{\delta}$, which is $1$-Lipschitz and within $\delta$-distance to $v$, thanks to the McShane formula, see \cite[Proposition 3.8]{Ciosmak20213}.

Observe that $A^c=\bigcup\{K_{\delta}^{\epsilon}\mid \delta>0\}\cup B$, where $B$ is the set of points belonging to at least two distinct transport rays of $v$, which is of Lebesgue measure zero.

Let us note that 
\begin{equation*}
    \frac{v-v_{\delta}}{\delta}\in\mathcal{F}
\end{equation*}
and that, thanks to the dominated convergence theorem,
\begin{equation*}
    \lim_{\epsilon\to 0}\lim_{\delta\to 0} \frac{v-v_{\delta}}{\delta}=\mathbf{1}_A\text{ in }L^1(\Omega,\mu+\nu).
\end{equation*}
This is to say, $\mathbf{1}_A\in\mathrm{cl}\mathcal{F}$. Applying an analogous reasoning, setting 
$v_{\delta,0}=v+\delta$ on $K^{\epsilon}$ and $v_{\delta,0}=v$ on $K_{\delta}^{\epsilon}$, shows that $-\mathbf{1}_A\in\mathrm{cl}\mathcal{F}$.
This is to say $\mathbf{1}_A\in\mathcal{A}$.

This shows that that the subspace of Borel measurable, integrable functions that are constant on transport rays of $v$ is contained in $\mathcal{A}$.

Let  us show that $\mathrm{cl}(\mathcal{F}\cap(-\mathcal{F}))$ is contained in the space of all Borel measurable, integrable functions that are constant on transport rays of $v$. 

If suffices to show that $\mathcal{F}\cap (-\mathcal{F})$ is contained in that space.

Let $u_1\colon \Omega\to\mathbb{R}$ be a $1$-Lipschitz function on $\Omega $ and let $\lambda_1> 0$ be such that 
\begin{equation*}
\lambda_1(v-u_1)\in\mathcal{F}\cap (-\mathcal{F}).
\end{equation*}
This holds if and only if there exists $1$-Lipschitz $u_2\colon\Omega\to\mathbb{R}$ and $\lambda_2> 0$ such that 
\begin{equation*}
    \lambda_1(v-u_1)=-\lambda_2(v-u_2).
\end{equation*}
Equivalently, there exists $\epsilon_0\in (0,1)$ such that \footnote{Note that we look at $u_1$ that lies in the Gleason part of $v$ of the set of $1$-Lipschitz functions.}
\begin{equation*}
    \frac{v-\epsilon_0  u_1}{1-\epsilon_0}\text{ is }1\text{-Lipschitz}.
\end{equation*}
Clearly, thanks to convexity of the  set of $1$-Lipschitz functions, this is equivalent to 
\begin{equation*}
    \frac{v-\epsilon  u_1}{1-\epsilon}\text{ is }1\text{-Lipschitz for all }\epsilon\in (0,\epsilon_0].
\end{equation*}
Note that, thanks to continuity, Lipschitzness needs to be verified merely on a set of full Lebesgue measure. Since the set  of points that belong to at least two distinct transport rays is of Lebesgue measure zero, see \cite[Corollary 2.15., p. 224]{Ciosmak2021}, it suffices to show that the inequality holds for the set $B^c$ of points that belong to a unique transport  ray. 
That is, we need that for all $\omega_1,\omega_2\in B^c$
\begin{equation*}
\abs{v(\omega_1)-v(\omega_2)-\epsilon(u_1(\omega_1)-u_1(\omega_2))}\leq (1-\epsilon)d(\omega_1,\omega_2).
\end{equation*}
Taking squares and rearranging we see that the desired inequality states that
\begin{equation*}
\epsilon^2\big(\abs{u_1(\omega_1)-u_1(\omega_2)}^2-d(\omega_1,\omega_2)^2\big)
+2\epsilon \big(d(\omega_1,\omega_2)^2-(v(\omega_1)-v(\omega_2))(u_1(\omega_1)-u_1(\omega_2))\big)
\end{equation*}
has to be bounded from the above by
\begin{equation*}
d(\omega_1,\omega_2)^2-\abs{v(\omega_1)-v(\omega_2)}^2.
\end{equation*}
In particular,  whenever $\abs{v(\omega_1)-v(\omega_2)}=d(\omega_1,\omega_2)$, then for all $0<\epsilon<\epsilon_0$ we see that  
\begin{equation*}
\epsilon\big(\abs{u_1(\omega_1)-u_1(\omega_2)}^2-d(\omega_1,\omega_2)^2\big)
+2 \big(d(\omega_1,\omega_2)^2-(v(\omega_1)-v(\omega_2))(u_1(\omega_1)-u_1(\omega_2))\big)
\end{equation*}
is non-positive. 
Letting $\epsilon$ tend to zero, we see that 
\begin{equation*}
d(\omega_1,\omega_2)^2=(v(\omega_1)-v(\omega_2))(u_1(\omega_1)-u_1(\omega_2)),
\end{equation*}
which in turn shows that 
\begin{equation*}
    u_1(\omega_1)-u_1(\omega_2)=v(\omega_1)-v(\omega_2)=d(\omega_1,\omega_2).
\end{equation*}
We see therefore that the transport rays of $v$ are contained in transport rays of $u_1$.
In particular, the function
\begin{equation*}
\lambda_1(v- u_1)
\end{equation*}
is constant on transport rays of $v$ and so  any limit of functions in $\mathcal{F}\cap(-\mathcal{F})$.

The proof is complete.
\end{proof}

\begin{remark}
    A version of the above proposition would also be available for pointwise closure, in place of the closure in $L^1(\mu+\nu)$, and with functions that are constant on the relative interia of transport rays, in place of functions constant on transport rays.
\end{remark}

If we then look at two measures $\mu,\nu$ for which $v$ is the optimal potential, then $\mu\prec_{\mathrm{cl}\mathcal{F}}\nu$. Let $B$ be the set of points that belongs to at least two distinct transport rays of $v$. It is a Borel set of Lebesgue measure zero, cf. \cite{Ciosmak2021}, thus we may consider the space $\Omega\setminus B$, in which the collection of all transport rays is pairwise disjoint.  This already  tells us that that transport rays form a partitioning of $\Omega$, thus the current discussion will be formal. Its purpose is to show that the localisation obtained via optimal transport can be seen as an instance of the general localisation scheme.

We shall apply Theorem \ref{thm:apir}. 
Let $Z\subset\mathcal{A}$ be a finite subset. Then the preimages $\Phi^{-1}(\mathrm{eirc}_{\mathcal{A}_Z}(\mu,\nu)(\cdot))$ of the irreducible components consist of transport sets, which are the fibres of elements of $Z$.
The partitioning into transport rays is countably Lipschitz, cf. \cite[Lemma 3.5., p. 226]{Ciosmak2021} for the Euclidean case and \cite[Proposition 3.2.1.]{Klartag2017}. Therefore picking appropriately the countable family of functions in $\mathcal{A}$ we easily may arrange that the intersections of such preimages give rise to a collection of all transport rays. Theorem \ref{thm:apir} shows now that these are pairwise disjoint.
\bibliographystyle{amsplain}
\bibliography{references}

\end{document}